\documentclass[12pt]{article}
\usepackage{amsthm,amsmath,amssymb,amsopn}
\usepackage{float}
\usepackage{booktabs}
\makeatletter
\def\@setparskip#1{%
	\setlength{\parskip}{#1}%
	\setlength{\topsep}{#1}%
	\setlength{\partopsep}{0pt}%
}
\@setparskip{0pt plus 0pt minus 1pt}
\makeatother
\usepackage{graphicx,epstopdf,epsfig,subfig}
\usepackage[colorlinks,linktocpage,linkcolor=black]{hyperref}
\usepackage{url,color}
\usepackage{booktabs}
\usepackage{colortbl,multirow}
\usepackage{appendix}

\numberwithin{equation}{section}
\numberwithin{table}{section}
\numberwithin{figure}{section}

\allowdisplaybreaks

\newtheorem{lemma}{Lemma}[section]
\newtheorem{theorem}{Theorem}[section]
\newtheorem{definition}{Definition}[section]
\newtheorem{remark}{Remark}[section]
\newtheorem{assumption}{Assumption}[section]

\newtheorem{claim}{Claim}[section]
\theoremstyle{definition}

\newcommand\bm{\boldsymbol}

\newcommand{\jump}[1]{[\![#1]\!]}

\newcommand{\lnorm}[1]{\Vert #1\Vert_{1,h}}
\newcommand{\sipnorm}[1]{\Vert #1\Vert_\sip}

\newcommand{\ddx}[2]{D_{(#1,#2)}}
\newcommand{\ccoef}[2]{c_{#1,#2}}
\newcommand{\gcoef}[2]{\gamma_{#1,#2}}

\def\dx{\,\mathrm{d}x} 

\def\dy{\,\mathrm{d}y}
\def\ds{\mathrm{d}s}

\def\sip{\mathrm{sip}}
\def\sips{{\mathrm{sip}{,}\ast}}

\def\bu{\bar{u}}

\def\bi{a^\sip_h}

\def\Rone{\uppercase\expandafter{\romannumeral1} }
\def\Rtwo{\uppercase\expandafter{\romannumeral2} }
\def\Rthree{\uppercase\expandafter{\romannumeral3} }
\def\Rfour{\uppercase\expandafter{\romannumeral4} }
\def\Rfive{\uppercase\expandafter{\romannumeral5} }
\def\Rsix{\uppercase\expandafter{\romannumeral6}}

\providecommand{\aave}[1]{ \{\mspace{-6.0mu}\{ #1 \}\mspace{-6.0mu}\} }
\providecommand{\jjump}[1]{ [\mspace{-2.5mu}[ #1 ]\mspace{-2.5mu}] }

\usepackage{mathtools}
\usepackage{stmaryrd}
\usepackage{authblk}
\usepackage{bm}

\makeatletter

\newcommand*\rel@kern[1]{\kern#1\dimexpr\macc@kerna}
\newcommand*\widebar[1]{%
	\begingroup
	\def\mathaccent##1##2{%
		\rel@kern{0.8}%
		\overline{\rel@kern{-0.8}\macc@nucleus\rel@kern{0.2}}%
		\rel@kern{-0.2}%
	}%
	\macc@depth\@ne
	\let\math@bgroup\@empty \let\math@egroup\macc@set@skewchar
	\mathsurround\z@ \frozen@everymath{\mathgroup\macc@group\relax}%
	\macc@set@skewchar\relax
	\let\mathaccentV\macc@nested@a
	\macc@nested@a\relax111{#1}%
	\endgroup
}
\makeatother

\begin{document}
\allowdisplaybreaks
\baselineskip=2pc

\begin{center}
{\bf Finite Volume Analysis of the Poisson Problem via a Reduced Discontinuous Galerkin Space}
\end{center}

\centerline{
Wenbo Hu\footnote{School of Mathematical Sciences, University of Science and Technology of China, Hefei, Anhui 230026, China. E-mail: huwenb@mail.ustc.edu.cn.},\;
Yinhua Xia\footnote{School of Mathematical Sciences, University of Science and Technology of China, Hefei, Anhui 230026, China. 
E-mail: yhxia@ustc.edu.cn. Research is supported in part by Anhui Provincial Natural Science Foundation grant No. 2408085J004, and NSFC grant No. 12271498}
}

\vspace{.2in}

\baselineskip=1.7pc
\centerline{\bf Abstract}

In this paper, we analyze the finite volume method for the Poisson equation by recasting it as a discontinuous Petrov-Galerkin scheme built upon the reduced discontinuous Galerkin (RDG) space. The key idea is to employ the RDG space as the trial space and piecewise constants as the test space, thereby formulating the method within a Petrov-Galerkin framework. This construction preserves the local conservation property characteristic of finite volume methods, while leveraging the approximation power of discontinuous Galerkin spaces at substantially reduced computational cost. An interesting observation is that the penalty parameter arising in this formulation differs markedly from its counterpart in standard discontinuous Galerkin methods. We provide a rigorous error analysis of the proposed scheme. In particular, we establish optimal-order convergence in the DG energy norm and suboptimal-order convergence rate in the \(L^2\) norm. The theoretical analysis is supported by a set of one- and two-dimensional numerical experiments with Dirichlet and periodic boundary conditions.  
This study highlights the role of the RDG space as an effective bridge between finite volume and discontinuous Galerkin methodologies, enabling finite volume schemes equipped with a mathematically rigorous convergence theory.

\vspace{.3cm}

\baselineskip=2pc

\section{Introduction}\label{sec:introduction}

The finite volume method (FVM) has long been a cornerstone in computational fluid dynamics and heat transfer simulations due to its inherent local mass conservation and geometric flexibility~\cite{godunov1959fvm,patankar1972fvm,patankar1980fvm,eymard2000fvm,vanLeer1979fvm,harten1987uniformly,liu1994weno,onate1993derivation,hampshire1992three}. Despite its widespread engineering success, establishing a rigorous theoretical framework for FVM has historically presented significant challenges. Unlike the finite element method (FEM), which rests on a clear variational foundation, the classical FVM was originally derived directly from integral conservation laws. Consequently, early error analyses were often dependent on geometric constraints, such as requiring admissible meshes where the line connecting cell centers is orthogonal to the interface. To overcome these limitations and establish a sound mathematical framework, theoretical developments have generally evolved along two distinct paths: vertex-centered and cell-centered formulations.

Research on vertex-centered schemes which are often referred to as finite volume element methods, is relatively mature, largely benefiting from their close relationship with linear finite element methods.
Seminal works in the 1980s and 1990s laid the foundation. Bank and Rose~\cite{bank1987fvm} established the rigorous theoretical framework for the vertex-centered box method applied to elliptic problems by interpreting it as a Petrov-Galerkin scheme. They demonstrated that, for the Poisson equation, the stiffness matrix is identical to that of linear FEM, thereby guaranteeing optimal convergence in the energy norm, comparable to that of standard linear finite elements. Subsequently, Cai~\cite{cai1991fvem} systematized this theory for general self-adjoint elliptic equations, deriving $O(h)$ error estimates in the discrete $H^1$ norm on shape-regular unstructured triangulations. Regarding the geometric constraints on triangular meshes, Mishev~\cite{mishev1998fvm} demonstrated that the classic vertex-centered scheme for elliptic problems yields optimal \(O(h)\) convergence and preserves monotonicity, provided that the underlying mesh satisfies the Voronoi property.
More recently, efforts have been extended to high-order formulations, where significant progress has been made through different theoretical approaches. Ye~\cite{ye2004fvm} proposed a discontinuous finite volume (DFV) method for elliptic problems, using discontinuous piecewise polynomials and penalty terms to ensure stability on minimal-support dual partitions, thereby establishing optimal \(O(h)\) convergence in the discrete $H^1$ norm. Building on this, Chou and Ye~\cite{ye2007fvm} developed a unified analysis framework covering conforming, nonconforming, and discontinuous finite volume schemes. Through a rigorous duality argument within this framework, they successfully derived the previously elusive optimal \(O(h^2)\) $L^2$ norm error estimates.
Building on the algebraic connections established by Bank and Rose, Xu and Zou~\cite{xu2009fvm} derived a generalized identity between the stiffness matrices of FVM and FEM. They established the inf-sup condition for both linear and quadratic vertex-centered schemes on simplicial meshes, thereby extending the optimal convergence theory to higher-order approximations. Chen~\cite{chen2010fvm} simplified the geometric construction by introducing a hierarchical decomposition strategy, applying finite volume constraints to vertex unknowns and Galerkin projections to high-order bubble functions. This hybrid structure ensures the inf-sup condition on unstructured meshes without requiring complex dual partitions for non-vertex degrees of freedom. Cao et al.~\cite{cao2013fvm} focused on enhancing accuracy by exploiting superconvergence. By defining trial functions on Gauss-Lobatto nodes and constructing dual partition boundaries at Gauss points, they achieved superior convergence rates at nodal points compared to standard global estimates.

Conversely, the theoretical landscape for cell-centered formulations, which constitutes the focus of this work, has evolved through a fundamentally different path, primarily relying on connections to the finite difference method (FD), the equivalence with mixed finite element methods and Petrov-Galerkin variational formulations. Forsyth et al.~\cite{forsyth1988fvm} utilized maximum principles to establish quadratic convergence for cell-centered schemes on smooth grids, effectively resolving earlier concerns regarding local truncation errors.  Jones and Menzies~\cite{jones2000fvm} addressed the inconsistency on nonuniform rectangular grids. They demonstrated that, although the local truncation error degrades to $O(1)$ due to mesh irregularity, the stability of the equivalent finite difference operator ensures \(O(h^2)\) global convergence. S\"uli~\cite{suli1991poisson} formulated the scheme as a Petrov-Galerkin method, establishing optimal \(O(h)\) convergence in the discrete $H^1$ norm on nonuniform meshes, while also demonstrating \(O(h^2)\)  convergence under higher regularity assumptions and mild mesh constraints. 
Parallel to this, Chavent and Roberts~\cite{chavent1991fvm} and Arbogast et al.~\cite{arbogast1997fvm} demonstrated that, specifically on rectangular grids, cell-centered finite volume schemes are algebraically equivalent to the lowest-order Raviart-Thomas ($RT_0$) mixed method with specific quadrature rules. Subsequently, Arbogast et al.~\cite{arbogast1998fvm} significantly extended this framework to geometrically general domains. They derived cell-centered schemes for full tensor coefficients on irregular meshes by utilizing an expanded mixed formulation. Although they introduced an enhanced technique with interface Lagrange multipliers to handle mesh distortions, their analysis remained confined to the lowest-order approximation.
Extending this variational equivalence to high-order schemes on general unstructured triangulations remains open.
On rectangular grids, however, substantial progress has been made. Cai et al.~\cite{cai2003fvm} systematically derived high-order finite volume schemes based on \(RTN\) and $BDFM$ mixed spaces. Crucially, they acknowledged the intrinsic limitations of this approach. For high-order \(BDFM\) spaces, a satisfactory finite volume reduction is currently unachievable. Furthermore, they noted that as the polynomial degree increases, the computational efficiency advantage of the derived FVM diminishes, observing that the resulting finite volume scheme may actually require more degrees of freedom than the underlying hybridized mixed finite element method.

Subsequently, theoretical research diverged. While unified frameworks like the gradient discretisation method \cite{droniou2018gdm} solidified the analysis for low-order schemes on general meshes, the pursuit of high-order accuracy increasingly converged towards DG methods. However, despite their rigorous variational foundation, DG methods typically incur high computational costs due to the large number of coupled degrees of freedom.

In this paper, we address this challenge by formulating the cell-centered finite volume scheme as a variational Petrov-Galerkin method to analyze its convergence properties for the Poisson model problem on rectangular grids. Specifically, we utilize the reduced discontinuous Galerkin (RDG) space introduced in \cite{hou2024rdg} as the trial space, while retaining the piecewise constant space as the test space. This formulation constitutes a high-order FVM. We rigorously prove that, on uniform meshes, this variational formulation is mathematically equivalent to the classical cell-centered FVM. This equivalence allows us to leverage the theories of DG methods to derive optimal error estimates in the DG energy norm and suboptimal error estimates in the \(L^2\) norm. Validated by numerical tests, this work provides a rigorous theoretical foundation for high-order FVMs on structured grids.

The rest of this paper is organized as follows. In Section~\ref{sec:notation}, we introduce the notations, definitions, and preliminaries. Section~\ref{sec:reconstructionspace} recalls the RDG space and its approximation properties. Section~\ref{sec:error_estimate} presents the proposed finite volume scheme for the Poisson problem, together with rigorous error estimates in the DG energy norm and \(L^2\) norm. Section~\ref{sec:numericalresults} reports numerical experiments that validate the theoretical findings. Finally, Section~\ref{sec:conclusion} concludes the paper with remarks and perspectives.

\section{Notation and Preliminaries}\label{sec:notation}

In this section, we introduce the notation, basic definitions of function spaces, and a set of auxiliary inequalities that will be used in the subsequent analysis. For clarity of exposition, we begin with the case of a hypercubic domain discretized by a Cartesian product mesh. We emphasize, however, that the reduced discontinuous Galerkin (RDG) space is not limited to this setting. Extensions to general domains and unstructured meshes will be addressed in future work. Throughout the paper, $C$ denotes a generic positive constant, independent of the mesh size $h$, whose value may change from line to line.

\subsection{Function spaces}
Let $K \subset \mathbb{R}^d$, $d=1,2$, be an open and bounded domain. For $s \geq 0$ and $1 \leq p \leq \infty$, we denote by $W^{s,p}(K)$ the Sobolev space of functions whose derivatives up to order $s$ belong to $L^p(K)$. The associated seminorm and norm are written as $|\cdot|_{s,p,K}$ and $\|\cdot\|_{s,p,K}$, respectively.  

In the special case $p=2$, we write $H^s(K):=W^{s,2}(K)$, with seminorm $|\cdot|_{s,K}$ and norm $\|\cdot\|_{s,K}$. For $s=0$, $H^0(K)$ coincides with $L^2(K)$, and we denote the norm and inner product by $\|\cdot\|_{0,K}$ and $(\cdot,\cdot)_K$, respectively.  

For $k \geq 0$, let $\mathcal{Q}^k(K)$ denote the space of tensor-product polynomials of degree at most $k$ on $K$.

\subsection{Meshes and notation}
Consider the $d$-dimensional hypercube $\Omega=\prod_{i=1}^d [0,1]$. Let $\mathcal{T}_h$ denote its uniform partition into $\mathcal{N}$ non-overlapping elements. In two dimensions, we denote the boundary segments by
\begin{align*}
\Gamma_{x=0}=\{(x,y)\in \partial\Omega:\,x=0\},\quad \Gamma_{x=1}=\{(x,y)\in \partial\Omega:\,x=1\},\\
\Gamma_{y=0}=\{(x,y)\in \partial\Omega:\,y=0\},\quad \Gamma_{y=1}=\{(x,y)\in \partial\Omega:\,y=1\}.
\end{align*}

Suppose $N_i$ elements are distributed along the $i$-th coordinate direction, i.e.
\[
[0,1]=\bigcup_{j=1}^{N_i} I_j,\qquad I_j=[x_{j-1},x_j].
\]
Then the total number of elements is $\mathcal{N}=\prod_{i=1}^d N_i$. For convenience, define the multi-index set
\[
\Theta^d=\{\bm{\theta}=(\theta_1,\dots,\theta_d):\, 1 \leq \theta_i \leq N_i, \; i=1,\dots,d\}.
\]

For each $\bm{\theta}\in\Theta^d$, the corresponding element centroid is denoted by
\begin{align}\label{eq:centroid_notation}
\bm{x}_{\bm{\theta}}=\big(x_{\bm\theta+1/2}^1,\dots,x_{\bm\theta+1/2}^d\big),
\end{align}
where each component $x_{\theta_i+1/2}^i$ is the midpoint of the $i$-th interval:
\[
x_{\bm\theta+1/2}^i = \tfrac{1}{2}(x_{\theta_i}^i + x_{\theta_i+1}^i).
\]

Let $\mathcal{F}_h$ denote the set of all element faces (edges when $d=2$), and let $\mathcal{F}_h^0 = \mathcal{F}_h \setminus \partial\Omega$ be the set of interior faces.  
For each $K \in \mathcal{T}_h$, we denote its measure by $|K|$ and its mesh size by $h_K$.  
In what follows, we restrict our attention to a uniform Cartesian mesh with mesh size $h$ and $N_i = N$ for simplicity.

To enforce boundary conditions, we introduce a layer of ghost (or virtual) cells surrounding the physical domain $\Omega$. The resulting extended mesh is defined by
\[
\mathcal{T}_h^{\text{ext}} \;=\; \mathcal{T}_h \cup \mathcal{T}_h^{\text{virt}}, 
\qquad 
\mathcal{T}_h^{\text{virt}} \;=\; \prod_{i=1}^d [-h,1+h] \setminus \mathcal{T}_h.
\]
The values of the numerical solution in $\mathcal{T}_h^{\text{virt}}$ are prescribed consistently with the boundary conditions. The specific treatment will be detailed in the subsequent sections for the different boundary types under consideration.

\subsection{Cell averages and projection operators}
For $u\in L^1(\mathcal{T}_h)$, its cell average on $K\in\mathcal{T}_h$ is defined as
\[
\bar{u}_K = \frac{1}{|K|}\int_K u\, \dx.
\]
\noindent
Specifically, for an element $K \in \mathcal{T}_h$ associated with the multi-index $\bm{\theta} \in \Theta^d$, we denote its cell average by its corresponding centroid index from \eqref{eq:centroid_notation}:
\[
 \bar{u}_K:=\bar{u}_{(\theta_1+1/2,\dots,\theta_d+1/2)} .
\]

This induces a projection operator $\Pi:L^1(\Omega)\to V^0_h$, where $V^0_h$ is the piecewise constant space, such that
\[
(\Pi u)|_K = \bar{u}_K,\qquad \forall K\in\mathcal{T}_h.
\]

We also define a discrete difference operator for cell averages. For $\bm{\theta}\in\Theta^d$ and $i=1,\dots,d$,
\[
D_{(\theta_1+1/2,\dots,\theta_i,\dots,\theta_d+1/2)} u =
\bar{u}_{(\theta_1+1/2,\dots,\theta_i+1/2,\dots,\theta_d+1/2)}-
\bar{u}_{(\theta_1+1/2,\dots,\theta_i-1/2,\dots,\theta_d+1/2)}.
\]

Using this, a mesh-dependent seminorm on $\mathcal{T}_h^{\text{ext}}$ is defined by
\[
|u|_{1,h}=\Big(\sum_{\bm{\theta}\in\Theta^d}\sum_{i=1}^d \big(D_{(\theta_1+1/2,\dots,\theta_i,\dots,\theta_d+1/2)}u\big)^2\Big)^{1/2}h^{d/2-1}.
\]

\subsection{Discrete function spaces}
On the mesh $\mathcal{T}_h$, the discontinuous Galerkin space of degree $k$ is
\[
V^k_h = \{v\in L^2(\Omega): v|_K\in \mathcal{Q}^k,\ \forall K\in\mathcal{T}_h\}.
\]
A subspace of $V^k_h$ is defined by
\[
\bar{V}^k_h=\{v\in V^k_h:\, \|v\|_{L^1(\Omega)}=0\}.
\]

We also employ broken Sobolev spaces:
\begin{align*}
H^s(\mathcal{T}_h)=\{v\in L^2(\Omega): v|_K\in H^s(K),\ \forall K\in\mathcal{T}_h\},\\
W^{s,p}(\mathcal{T}_h)=\{v\in L^2(\Omega): v|_K\in W^{s,p}(K),\ \forall K\in\mathcal{T}_h\}.
\end{align*}

For $u\in H^s(\mathcal{T}_h)$ or $W^{s,p}(\mathcal{T}_h)$, the broken gradient and Laplacian are defined elementwise:
\[
\nabla_h u := \nabla(u|_K),\quad \Delta_h u := \Delta(u|_K),\qquad \forall K\in\mathcal{T}_h.
\]

\subsection{Interface operators and norms}
For any interior face $F\in \mathcal{F}_h^0$ shared by elements $T_1$ and $T_2$, and for a function $v$ with well-defined traces, the average and jump are defined as
\[
\aave{v}_F(x)=\tfrac{1}{2}\big(v|_{T_1}(x)+v|_{T_2}(x)\big),\qquad 
\llbracket v\rrbracket_F(x)=v|_{T_1}(x)-v|_{T_2}(x).
\]
For vector-valued functions, these operators are applied componentwise. When the context is clear, the subscript $F$ and the variable $x$ will be omitted.

We define the symmetric interior penalty (sip) norm by
\[
\|v\|_{\text{sip}}^2=\sum_{K\in\mathcal{T}_h}\|\nabla v\|_{L^2(K)}^2+\sum_{F\in\mathcal{F}_h}\frac{1}{h}\|\llbracket v\rrbracket\|_{L^2(F)}^2,
\]
and the extended $\text{sip},*$ norm by
\[
\|v\|_{\text{sip},*}^2=\|v\|_{\text{sip}}^2+\sum_{K\in\mathcal{T}_h}\|\nabla v|_{K}\!\cdot\! \mathrm n_K\|_{L^2(\partial K)}^2,
\]
where $n_K$ is the unit outward normal on $\partial K$.

\subsection{Auxiliary inequalities}
We recall several standard inequalities used later in the analysis.

\begin{itemize}
\item \textbf{Young's inequality.} For $a,b\geq 0$ and $p,q>1$ with $1/p+1/q=1$,
\[
ab \leq \tfrac{a^p}{p}+\tfrac{b^q}{q}.
\]

\item \textbf{Agmon's inequality.} For any $v\in H^1(K)$,
\[
\|v\|_{L^2(\partial K)}^2 \leq C\big(h_K^{-1}\|v\|_{L^2(K)}^2 + h_K\|\nabla v\|_{L^2(K)}^2\big).
\]

\item \textbf{Approximation property.} For any $v\in H^{k+1}(K)$, there exists $v_h\in \mathcal{Q}^k(K)$ such that
\[
\|v-v_h\|_{0,K}+h|v-v_h|_{1,K}\leq C h^{k+1}|v|_{k+1,K}.
\]

\item \textbf{Inverse inequality.} For any $v_h\in V^k_h$,
\[
\|\partial_x v_h\|_{0,K}\leq C h^{-1}\|v_h\|_{0,K},\qquad \forall K\in\mathcal{T}_h.
\]
\end{itemize}

\section{Reduced Discontinuous Galerkin Space}\label{sec:reconstructionspace}

In this section, we introduce the reduced discontinuous Galerkin (RDG) operator and the corresponding RDG space, following the framework proposed in \cite{hou2024rdg}. We recall the well-posedness condition and summarize key approximation properties in the $L^2$ and $H^1$ norms. In addition, we establish an error estimate in the $\sips$ norm, which plays a crucial role in the error analysis of the proposed scheme.

\subsection{Reconstruction operator}

We begin by recalling the construction of the local stencil. For each element $K \in \mathcal{T}_h$, a stencil $S(K)$ is defined, consisting of $K$ together with its neighboring elements, chosen by standard strategies such as Moore-neighborhood selection. We denote
\[
S(K) = \{K_{\theta_1},K_{\theta_2},\dots,K_{\theta_s}\}, \qquad d(S(K)) = s,
\]
where $s$ is the number of elements in the stencil.

Let $\{L_i(x)\}_{i=0}^k$ be the one-dimensional Legendre basis functions of degree $i$ for $[-1,1]$. The tensor-product multidimensional basis on $K$ is
\[
L^\alpha_K(\bm x) = \prod_{i=1}^d \widehat{L}^{\alpha_i}\!\left(\frac{2(x_i-x_{\theta_i})}{h^i_K}\right),
\]
where $\bm\alpha=(\alpha_1,\dots,\alpha_d)\in A^k:=\{0\leq \alpha_i \leq k\}$, and $h^i_K$ is the element size in the $i$-th direction.  

For a function $u\in C^0(K)$, the Legendre moments are defined as
\[
I_K^{\bm\alpha} u(\bm x) = \Big(\prod_{i=1}^d \tfrac{2\alpha_i+1}{h^i_K}\Big) \int_K L_K^{\bm\alpha}(\bm x)\, u(\bm x)\, \mathrm{d}\bm x.
\]
These moments are used to construct a polynomial $p \in \mathcal{Q}^k(K)$ via the interpolation conditions
\[
I^{\bm\alpha}_K p(\bm x) = I^{\bm\alpha}_K u(\bm x), \qquad \forall {\bm\alpha}\in A^k.
\]
This defines the local $L^2$-projection operator $\pi^k_K : C^0(K) \to \mathcal{Q}^k(K)$.

The local reconstruction operator $R^k_K$ builds a polynomial of degree $k$ from lower-order approximations $\pi_{K_{\bm\theta}}^m u$ ($m<k$) on neighboring cells $K_{\bm\theta} \in S(K)$. It is defined by
\[
I^{\bm\alpha}_{K_{\bm\theta}} (R^k_K u) = I^{\bm\alpha}_{K_{\bm\theta}} u, \qquad \forall K_{\bm\theta} \in S(K), \ \alpha\in A^m,\ m<k.
\]
The global reconstruction operator $R^k$ is then defined elementwise by
\[
(R^k u)|_K = R^k_K u, \qquad \forall K\in \mathcal{T}_h.
\]

\subsection{Properties of the reconstruction operator}

\begin{assumption}\label{assumption1}
For every $K\in \mathcal{T}_h$ and every $w\in \mathcal{Q}^k(S(K))$, if the Legendre moments vanish, i.e.
\[
I^{\bm\alpha}_{K_{\bm\theta}}(R^k_K w)=0,\qquad \forall K_{\bm\theta} \in S(K),\ {\bm\alpha}\in A^m,
\]
then $w|_{S(K)}\equiv 0$.
\end{assumption}

This assumption guarantees uniqueness of the reconstruction and will be assumed henceforth.

\begin{lemma}[{\cite[Theorem 1]{hou2024rdg}}]
If Assumption \ref{assumption1} holds, then the reconstruction is $k$-exact:
\[
R^k_K u = u, \qquad \forall u\in \mathcal{Q}^k(S(K)).
\]
Moreover, if $u\in H^{k+1}(S(K))$, then
\begin{align}
\|u-R^k_K u\|_{L^2(K)} &\leq C h^{k+1} \|u\|_{H^{k+1}(S(K))}, \label{eq:locall2}\\
|u-R^k_K u|_{H^1(K)} &\leq C h^k \|u\|_{H^{k+1}(S(K))}, \label{eq:localh1}
\end{align}
where $C$ is independent of $h$.
\end{lemma}

\begin{theorem}\label{th:interpolationestimate}
If Assumption \ref{assumption1} holds, then for $u\in H^{k+1}(\Omega)$ the global reconstruction $R^k u$ satisfies
\begin{align}
\|u-R^k u\|_{L^2(\Omega)} &\leq C h^{k+1}, \label{eq:globall2}\\
|u-R^k u|_{H^1(\Omega)} &\leq C h^k, \label{eq:globalh1}\\
\|u-R^k u\|_{\sips} &\leq C h^k, \label{eq:globalsips}
\end{align}
where $C$ is independent of $h$.
\end{theorem}

\begin{proof}
Estimate \eqref{eq:globall2} follows directly:
\[
\|u-R^k u\|_{L^2(\Omega)}^2 \leq C \sum_{K\in \mathcal{T}_h} \|u\|_{H^{k+1}(S(K))}^2 h^{2k+2}.
\]
Since the stencil size $\# S(K)$ is fixed, we have
\[
\sum_{K\in \mathcal{T}_h} \|u\|_{H^{k+1}(S(K))}^2 \leq C \|u\|_{H^{k+1}(\Omega)}^2,
\]
yielding \eqref{eq:globall2}. Estimate \eqref{eq:globalh1} is analogous. Finally, combining Agmon's inequality, the definition of the $\sips$ norm, and the interpolation bounds \eqref{eq:globall2}-\eqref{eq:globalh1} gives \eqref{eq:globalsips}.
\end{proof}

\subsection{Reduced discontinuous Galerkin space}

Let $V^m_h$ denote the standard DG space of piecewise polynomials of degree $m$ with $m<k$. The global reconstruction $R^k$ maps $V^m_h$ into $V^k_h$, yielding a reduced high-order space
\[
U^k_h := R^k V^m_h.
\]
This space is spanned by reconstructed basis functions
\[
\phi^{\bm\alpha}_K = R^k v^{\bm\alpha}_K, \qquad K\in\mathcal{T}_h,\ {\bm\alpha}\in A^m,
\]
where $v^{\bm\alpha}_K$ is a locally supported function associated with Legendre moment ${\bm\alpha}$. Thus
\[
U^k_h = \text{span}\{\phi^{\bm\alpha}_K : K\in\mathcal{T}_h,\ {\bm\alpha}\in A^m\}.
\]

The RDG space $U^k_h$ is a subspace of $V^k_h$, inheriting its approximation power but with significantly fewer degrees of freedom: $O(N \cdot \dim(A^m))$ versus $O(N \cdot \dim(A^k))$ for $V^k_h$. This reduction in dimensionality leads to improved computational efficiency without compromising high-order accuracy.

\section{Analysis of the Finite Volume Method for the Poisson Problem}\label{sec:error_estimate}

In this section, we demonstrate that the finite volume method for the Poisson problem can be interpreted as a special case of the IPDG formulation, where the trial and test spaces are chosen as $U_h^k$ and $V_h^0$, respectively. Based on this bilinear form, we then establish the corresponding error estimates. For clarity of presentation, we focus on the second-order case, while noting that the analysis can be extended to higher-order formulations. We first consider the homogeneous Dirichlet boundary condition and subsequently extend the results to the periodic problem.

\subsection{Dirichlet Problem}\label{sec:dirichlet_ipdg}

Consider the Poisson equation: find $u \in H^2(\Omega)$ such that
\begin{equation} \label{eq:1d_dirichlet_poisson}
	\begin{cases}
		-\Delta u = f, & x \in \Omega,\\
		u = 0, & x \in \partial\Omega,
	\end{cases}
\end{equation}
where $f \in L^2(\Omega)$ and $\Omega=[0,1]^d$.  We first show that the finite volume scheme for the Poisson problem can be written as the IPDG scheme with the help of the reduced DG space.

The symmetric IPDG bilinear form~\cite{arnold2002ipdg,daniele2012ipdg} is
\begin{equation}
\begin{aligned}\label{eq:bi1}
	a_h^{\sip}(u,v_h) &= -\sum_{K\in\mathcal{T}_h}\int_K \nabla_h u\!\cdot\!\nabla_h v_h \,\dx
	-\sum_{F\in\mathcal{F}_h}\int_F \aave{\nabla_h u}\!\cdot\! \mathbf{n}_F \,\jjump{v_h}\,\ds \\
	&\quad -\sum_{F\in\mathcal{F}_h}\int_F \jjump{u}\,\aave{\nabla_F v_h}\!\cdot\! \mathbf{n}_F \,\ds
	+\sum_{F\in\mathcal{F}_h}\frac{\eta}{h_F}\int_F \jjump{u}\jjump{v_h}\,\ds,
\end{aligned}
\end{equation}
or equivalently,
\begin{equation}
\begin{aligned}\label{eq:bi2}
	a_h^{\sip}(u,v_h) &= -\sum_{K\in\mathcal{T}_h}\int_K (\Delta u)v_h \,\dx
+\sum_{F\in\mathcal{F}^i_h}\int_F \jjump{\nabla_h u}\!\cdot\! \mathbf{n}_F \,\aave{v_h}\,\ds \\
&\quad -\sum_{F\in\mathcal{F}_h}\int_F \jjump{u}\,\aave{\nabla_F v_h}\!\cdot\! \mathbf{n}_F \,\ds
+\sum_{F\in\mathcal{F}_h}\frac{\eta}{h_F}\int_F \jjump{u}\jjump{v_h}\,\ds,
\end{aligned}
\end{equation}
with the penalty parameter $\eta \ge 0$, 
and
\[
l(v_h) = \sum_{K\in \mathcal{T}_h}\int_K f v_h\,\dx.
\]

We take $V_h^r$ as the test space and $U_h^k = R^k V_h^r$ as the trial space.  
The resulting IPDG scheme reads: Find $u_h \in U_h^2$ such that
\begin{equation}\label{eq:1d_dirichlet_disproblem}
	a_h^{\sip}(u_h,v_h) = l(v_h), \quad \forall v_h \in V_h^0.
\end{equation}
We will show that for $(r,k) = (0,2)$ this reduces to a second-order finite volume method. 
In $U_h^2$, the boundary condition is strongly enforced by prescribing the numerical solution in the virtual cells $\mathcal{T}_h^{\text{virt}}$.  
With this choice of trial and test spaces, it is worth noting that the discrete scheme remains well-posed even for $\eta = 0$, as will be shown later.  
We first illustrate this in the one-dimensional case, where boundary conditions are enforced in the reconstruction operator by prescribing the cell averages in the virtual cells.

\begin{definition}
In the one-dimensional case, the reconstruction operator \(R^2\) is naturally defined on interior elements.  
For boundary elements, e.g., \([x_0,x_1]\) or \([x_{N-1},x_N]\), the definition of the local reconstruction operator \(R^2_K\) requires the introduction of virtual elements.  

As an example, consider the left boundary cell \(K=[x_0,x_1]\).  
We define \(R^2_K u \in \mathcal{P}^2(K)\) such that
\[
\int_{K'} R^2_K u \, \dx \;=\; \int_{K'} u \, \dx,
\qquad \forall K'\subset S(K)=\{[x_{-1},x_0],[x_0,x_1],[x_1,x_2]\}, 
\quad K'\cap\Omega\neq\varnothing.
\]
Here, the cell average on the virtual cell \([x_{-1},x_0]\) is prescribed to enforce the boundary condition in \eqref{eq:1d_dirichlet_poisson}:
\[
R^2_K u(0) = 0 \quad\text{and}\quad R^2_K u(1) = 0. 
\]
\end{definition}

Next, we show that for functions in the trial space $U_h^2$, the internal jump of the discrete gradient vanishes.

\begin{lemma}\label{le:1d_grad_jump}
	In the one-dimensional case, for any $u \in V_h^0$, the reconstructed function $R^k u \in U_h^k$ satisfies
	\[
	\jjump{\nabla_h (R^k u)} = 0,
	\]
	provided that the reconstruction order $k$ is even.
\end{lemma}

\begin{proof}
	Consider the element $[x_i, x_{i+1}]$ mapped to the reference interval $[0,1]$.  
	On this reference element, the reconstructed function can be expressed as
	\[
	\widehat{(R^k u)}|_{[0,1]}(\hat{x})
	= \sum_{k=-m}^{m} \Big( \sum_{s=-m}^{k} \bar{u}_{i+s+1/2} \Big) L'_{k+1}(\hat{x}),
	\]
	where $L_k(\hat{x})$ denotes the Lagrange basis polynomial defined by
	\[
	L_k(\hat{x}) = \frac{\prod_{r=-m,\, r\neq k}^{m+1} (\hat{x}-r)}{\prod_{r=-m,\, r\neq k}^{m+1} (k-r)}.
	\]

	For $j=-m,\dots,m$, one can verify that
	\[
	\int_j^{j+1} \widehat{(R^k u)}|_{[0,1]}(\hat{x}) \, d\hat{x}
	= \bar{u}_{i+j+1/2},
	\]
	which satisfies the defining property of $R^k$.

	The derivative at the left interface $\hat{x}=0$ is
	\[
	\widehat{(R^k u)_{\hat{x}}}|_{[0,1]}(0)
	= \sum_{k=-m}^{m} \Big( \sum_{s=-m}^{k} \bar{u}_{i+s+1/2} \Big) L''_{k+1}(0)
	= \sum_{s=-m}^{m} \ccoef{s}{m} \, \bar{u}_{i+s+1/2},
	\]
	where
	\[
	\ccoef{s}{m} = \sum_{k=s}^{m} L''_{k+1}(0).
	\]
	Since $L_{m+1}(\hat{x})$ is an odd function, $L''_{m+1}(0)=0$, implying $\ccoef{m}{m}=0$ and hence
	\[
	\widehat{(R^k u)_{\hat{x}}}|_{[0,1]}(0)
	= \sum_{s=-m}^{m-1} \ccoef{s}{m} \, \bar{u}_{i+s+1/2}.
	\]

	Similarly, for the adjacent element $[x_{i-1}, x_i]$ mapped to $[0,1]$, we obtain
	\[
	\widehat{(R^k u)_{\hat{x}}}|_{[0,1]}(1^-)
	= \sum_{s=-m}^{m} \tilde c_{s,m} \, \bar{u}_{i+s-1/2},
	\qquad
\tilde c_{s,m} = \sum_{k=s}^{m} L''_{k+1}(1).
	\]
	Noting that $\tilde c_{-m,m} = 0$, we rewrite
	\[
	\widehat{(R^k u)_{\hat{x}}}|_{[0,1]}(1^-)
	= \sum_{s=-m}^{m-1} \tilde c_{s+1,m} \, \bar{u}_{i+s+1/2}.
	\]

	Finally, we show that $\ccoef{s}{m} = \tilde c_{s+1,m}$, for $s = -m, \dots, m-1$.
	Define
	\[
	S_{s+1}(x) = \sum_{j=s+1}^{m+1} L_j(x),
	\]
	which satisfies
	\[
	S_{s+1}(l) = 
	\begin{cases}
		0, & l = -m, \dots, s, \\
		1, & l = s+1, \dots, m+1.
	\end{cases}
	\]
	Let $T_s(x) = S_{s+1}(x-1)$, then $T_s''(1) = S_{s+1}''(0)$.
	Expressing $T_s$ in the Lagrange basis gives
	\[
	T_s(x) = \sum_{j=s+1}^{m} L_{j+1}(x),
	\]
	and thus
	\[
	\ccoef{s}{m} = S''_{s+1}(0) = T_s''(1) = {\tilde c_{s+1,m}}.
	\]
	This establishes the continuity of the reconstructed derivative, i.e.,
	\(
	\jjump{\nabla_h (R^k u)} = 0.
	\)
\end{proof}


This lemma implies that the second term in the bilinear form~\eqref{eq:bi2} vanishes for functions in the trial space $U_h^k$.  
Furthermore, when the test space is chosen as $V_h^0$, the third term also vanishes.  
Consequently, by setting $\eta = 0$ in the penalty term, the IPDG scheme \eqref{eq:1d_dirichlet_disproblem} reduces exactly to the classical second-order finite volume scheme for the Poisson problem. Now, we establish the core lemma that describes the coefficients of the reconstructed derivative flux and their properties, which is essential for the stability proof.
\begin{lemma}\label{lem:1d_deriv_representation} 
	For one dimensional case, there exist coefficients \(\left(\gcoef{k}{m}\right)_{k=-m+1}^{m-1}\), depending only on \(m\), such that for any \(u \in H^{2m+1}(\Omega)\), the reconstructed derivative at node \(x_i\) is given by
	\begin{equation}\label{eq:deriv_representation}
		(R^ku)_x(x_i) = \frac{1}{h} \sum_{k=-m+1}^{m-1} \gcoef{k}{m} D_{(i+k)}u.
	\end{equation}
	Furthermore, these coefficients satisfy the property:
	\begin{equation}\label{eq:coeff_property}
		\gcoef{0}{m} > \sum_{k=-m+1, k\neq 0}^{m-1} |\gcoef{k}{m}|.
	\end{equation}
\end{lemma}

We refer the reader to Appendix~\ref{ap:proof_of_le:1d_deriv_representation} for the detailed proof.

\subsubsection{1D Error Estimates}\label{subsec:1d_dirichlet_error}
In this section, we derive error estimates for the one-dimensional scheme \eqref{eq:1d_dirichlet_disproblem}.   With the boundary condition strongly enforced in the reconstruction operator, we obtain the following theorem establishing the equivalence of the discrete norms.

\begin{lemma}[Norm equivalence]\label{eq:1dnormequivalence}
    Let \(U^2_h\) be equipped with the semi-norm \(\vert \cdot \vert_{1,h}\) and the norm \(\|\cdot\|_\sip\). Then \(\vert \cdot \vert_{1,h}\) defines a norm on \(U^2_h\), which we denote by \(\|\cdot\|_{1,h}\). Furthermore, the norms \(\|\cdot\|_{1,h}\) and \(\|\cdot\|_\sip\) are equivalent, in the sense that there exist positive constants \(C_1\) and \(C_2\) such that
    \begin{align}
        C_1 \| u \|_{\mathrm{sip}} \leq \| u \|_{1,h} \leq C_2 \| u \|_{\mathrm{sip}}, \quad \forall u \in U^2_h.
    \end{align}
\end{lemma}

\begin{proof}
    First, we show that the semi-norm \(\vert \cdot \vert_{1,h}\) defines a norm on \(U^2_h\). The properties of non-negativity and homogeneity are readily satisfied. We now focus on establishing positivity.

    Suppose \(\| u \|_{1,h} = 0\). By definition, this implies that \(D_{(i)}u = 0\) for all \(i = 0, \ldots, N\). Consequently, we have \(\bar{u}_{i-1/2} \equiv C\) for all \(i = 0, \ldots, N+1\). The cell averages on the virtual elements are determined by the boundary conditions. Applying the Dirichlet boundary condition along with the definition of \(R^2\), we obtain the following equations:
    \begin{equation}\label{boundary_equality}
        \begin{cases}
        	\tfrac{1}{3}\bar{u}_{-1/2} + \tfrac{5}{6}\bar{u}_{1/2} - \tfrac{1}{6}\bar{u}_{3/2} = 0,\\ 
        	-\tfrac{1}{6}\bar{u}_{N-3/2} + \tfrac{5}{6}\bar{u}_{N-1/2} + \tfrac{1}{3}\bar{u}_{N+1/2} = 0.
        \end{cases}
    \end{equation}

    From this, we deduce \(\bar{u}_{i+1/2} = 0\) for all \(i = -1, \ldots, N\). Hence, the reconstructed function \(u = R^2 \bar{u} \equiv 0\). Therefore, we have proven that \(\vert \cdot \vert_{1,h}\) is indeed a norm on the space \(U^2_h\), which we denote as \(\|\cdot\|_{1,h}\).

    Next, we demonstrate the equivalence of \(\|\cdot\|_{1,h}\) and \(\|\cdot\|_{\sip}\). To establish this, we will derive bounds for each norm in terms of the other.

    For the first inequality, by the definition of \(\|\cdot\|_\sip\), we have:
    \begin{equation}
        \begin{aligned}
            \| u \|_{\mathrm{sip}}^2 & \leq \frac{1}{h}\sum_{j=0}^{N-1}\frac{1}{2}(D_{(j)}u^2 + D_{(j+1)}u^2) + \frac{1}{36h}\sum_{j=1}^{N-1}3(D_{(j-1)}u^2 + D_{(j)}u^2 + D_{(j+1)}u^2) \\
            & \leq \frac{C}{h}\sum_{j=0}^{N}D_{(j)}u^2 = C_1^2 \| u \|_{1,h}^2.
        \end{aligned}
    \end{equation}

    For the second inequality, using the definition of \(\| \cdot \|_{1,h}\) and the elementary inequality derived from \((D_{(j)}u - D_{(j+1)}u)^2 \geq 0\), we obtain the following estimate:
    \begin{equation}\label{eq:1d_dirichlet_1h<sip}
        \begin{aligned}
            \| u \|_{1,h}^2 & \leq \frac{1}{h}\sum_{j=0}^{N-1}(D_{(j)}u^2 + D_{(j+1)}u^2) \\
            & \leq \frac{1}{h}\sum_{j=0}^{N-1}\left((D_{(j)}u^2 + D_{(j+1)}u^2)+(D_{(j)}u^2 + D_{(j+1)}u^2 - 2 | D_{(j)}u D_{(j+1)}u |)\right) \\
            & \leq \frac{2}{h}\sum_{j=0}^{N-1}(D_{(j)}u^2 + D_{(j)}u D_{(j+1)}u + D_{(j+1)}u^2)\leq C_2^2 \| u \|_{\sip}^2.
        \end{aligned}
    \end{equation}

    Thus, we conclude that the norms \(\| u \|_{1,h}\) and \(\| u \|_\sip\) are indeed equivalent.
\end{proof}

\begin{lemma}[Reconstruction boundedness] \label{eq:1doperatorinjective}
    For all \( u\in V_h^0 \), the following inequality holds:
    \[
    \|u\|_\sip \leq C \|R^2 u\|_\sip,
    \]
    where \( R^2 u \in U^2_h \) and \( C \) is a positive constant independent of \( u \) and \( h \).
\end{lemma}

\begin{proof}
    By the definition of the \(\sip\) norm, for \( u \in V^0_h \), we have:
    \begin{align*}
        \|u\|_\sip^2 &= \sum_{i=1}^{N-1} \frac{1}{h} \jump{u}^2 = \sum_{i=1}^{N-1} \frac{1}{h} D_{(i)} u^2.
    \end{align*}
    
    By the properties of the reconstruction operator \( R^2 \), it follows that
    \[
    \widebar{u}_{i-1/2} = \widebar{R^2 u}_{i-1/2}, \quad \forall i=0, 1, \ldots, N+1.
    \]
    Because \( R^2 u \in U^2_h \), we can apply the bound given in \eqref{eq:1d_dirichlet_1h<sip}. Thus, we obtain:
    \[
    \|u\|_\sip^2 = \sum_{i=1}^{N-1} \frac{1}{h} (D_{(i)} u)^2 \leq C \| R^2 u \|_\sip^2,
    \]
    where \( C \) is a constant independent of \( u \) and \( h \). 

    Therefore, we conclude that 
    \[
    \| u \|_\sip \leq C \| R^2 u \|_\sip.
    \]
\end{proof}

With the definitions of the \(\sip\) norm and leveraging Lemmas \ref{eq:1dnormequivalence} and \ref{eq:1doperatorinjective}, we can demonstrate that the bilinear form \(\bi(\cdot,\cdot)\) satisfies both the inf--sup condition and the boundedness condition.

\begin{lemma}[Inf--sup condition]
	For the bilinear form \(a_h^\sip: U^2_h \times V^0_h \to \mathbb{R}\), the inf--sup condition holds:
	\begin{equation}
		\adjustlimits\inf_{u \in U_h^2 \setminus \{0\}} \sup_{v \in V_h^0 \setminus \{0\}} \frac{\vert a(u,v) \vert}{\| u \|_\sip \| v \|_\sip} \geq C,
		\label{eq:1d_dirichlet_infsupcondition}
	\end{equation}
	where \(C\) is a positive constant independent of \(\mathcal{T}_h\) and \(u, v\).
\end{lemma}

\begin{proof}
	We prove the condition for the specific case $(r,k) = (0,2)$, which, as noted in \eqref{eq:1d_dirichlet_disproblem}, defines our scheme. This choice corresponds to using the piecewise constant space $V_h^0$ as the test space and the 2nd-degree polynomial reconstruction space $U_h^2$ as the trial space. In this case, \(m=1\), and the summation indices in the definition of the bilinear form simplify significantly. 
		Let \(u \in U^2_h\) and \(v \in V^0_h\). By direct computation, we obtain the following expression for the bilinear form: 
	\begin{align*}
		a_h^\sip(u,v)&=\sum_{j=1}^{N-1}\left(\frac{1}{h}\sum_{k=-m+1}^{m-1} \gcoef{k}{m}D_{(j+k)}u\right) \, D_{(i+k)}v\\
		&\quad+\left(\frac{1}{h}\sum_{k=-m+1}^{m-1} \gcoef{k}{m}D_{(k)}u\right)\, v_{1/2}-\left(\frac{1}{h}\sum_{k=-m+1}^{m-1} \gcoef{k}{m}D_{(N+k)}u\right)\, v_{N-1/2} \\
		&\quad+\frac{\eta}{6h}\sum_{j=1}^{N-1}(-D_{(j-1)}u+2D_{(j)}u-D_{(j+1)}u)D_{(j)}v.
	\end{align*}
	Let \(v = \Pi u \in V^0_h\). Substituting the equality \eqref{boundary_equality} derived from the Dirichlet boundary condition into our expression yields:
	\begin{align*}
		a_h^\sip(u,v)&=\sum_{j=1}^{m-1}\left(\frac{1}{h}\sum_{k=-m+1}^{m-1} \gcoef{k}{m}D_{(j+k)}u\right) \, D_{(j+k)}u\\
		&\quad+\frac{1}{6}\left(\frac{1}{h}\sum_{k=-m+1}^{m-1} \gcoef{k}{m}D_{(k)}u\right)\, (2D_{(0)}u+D_{(1)}u)\\
		&\quad+\frac{1}{6}\left(\frac{1}{h}\sum_{k=-m+1}^{m-1} \gcoef{k}{m}D_{(N+k)}u\right)\, (D_{(N-1)}u+2D_{(N)}u) \\
	&\quad+\frac{\eta}{6h}\sum_{j=1}^{N-1}(-D_{(j-1)}u+2D_{(j)}u-D_{(j+1)}u)D_{(j)}u.
	\end{align*}
	Simplifying by grouping identical difference terms, we arrive at
		\begin{align*}
			a_h^{\sip}(u,\Pi u)
			&= \frac{1}{h}\left(
			\frac{1}{3}D_{(0)}u^2
			+\Big(1+\frac{\eta}{3}\Big)\sum_{j=1}^{N-1} D_{(j)}u^2
			+\frac{1}{3}D_{(N)}u^2\right.\\
			&\qquad\qquad
			\left.+\frac{1-\eta}{6}\,D_{(0)}u\,D_{(1)}u
			-\frac{\eta}{3}\sum_{j=1}^{N-2} D_{(j)}u\,D_{(j+1)}u
			+\frac{1-\eta}{6}\,D_{(N-1)}u\,D_{(N)}u	\right).
		\end{align*}
		
		By applying Young's inequality along with Lemmas~\ref{lem:1d_deriv_representation},~\ref{eq:1dnormequivalence} and~\ref{eq:1doperatorinjective}, we find that for all \(u \in U^2_h\):
\begin{align*}
	a_h^{\sip}(u,\Pi u)
	&\ge \frac{1}{h} \sum_{j=0}^{N} C_j D_{(j)}u^2 ,
\end{align*}
where
\[
C_j =
\begin{cases}
	\dfrac{1}{3} - \dfrac{|1 - \eta|}{12}, & j = 0,\, N, \\[1.2ex]
	\left(1 + \dfrac{\eta}{3}\right) - \dfrac{|1 - \eta|}{12} - \dfrac{|\eta|}{6}, & j = 1,\, N - 1, \\[1.2ex]
	\left(1 + \dfrac{\eta}{3}\right) - \dfrac{|\eta|}{3}, & 2 \le j \le N - 2.
\end{cases}
\]

By choosing \(\eta\in (-\frac{3}{2},5)\), we derive that
\begin{align*}
a_h^{\sip}(u,\Pi u)	&	\geq C \Vert u\Vert_{1,h}^2\geq C\Vert u\Vert_\sip\Vert R^2\Pi u\Vert_\sip\geq C\Vert u\Vert_\sip \Vert \Pi u\Vert_\sip.
\end{align*}
	
	Since \(v \in V^0_h\) is arbitrary, we conclude that
	\begin{align*}
	\adjustlimits\inf_{u \in U^2_h \setminus \{0\}} \sup_{v \in V^0_h \setminus \{0\}} \frac{\vert a(u,v) \vert}{\| u \|_\sip \| v \|_\sip} \geq C.
	\end{align*}
\end{proof}
\begin{remark}
	The above analysis ensures \(C_j > 0\) for \(-\tfrac{3}{2} < \eta < 5\), which guarantees the coercivity of \(a_h^{\sip}(\cdot,\cdot)\). 
	Nevertheless, numerical experiments indicate that this property may hold for a wider range of~\(\eta\).
\end{remark}

\begin{lemma} [Boundedness]\label{le:1d_dirichlet_boundedness}
    The bilinear form \(a^\sip_h(\cdot,\cdot) : (H^2 + U^2_h) \times V^0_h \to \mathbb{R}\) is bounded. Specifically, there exists a constant \(C > 0\), independent of \(h\), such that
    \begin{equation}
        \vert a^\sip_h(u,v) \vert \leq C \| u \|_\sips \| v \|_\sip, \quad \forall u \in U^2_h, \forall v \in V^0_h,
        \label{eq:1d_dirichlet_boundedness}
    \end{equation}
\end{lemma}

\begin{proof}
    By applying Lemma 4.16 from \cite{daniele2012ipdg}, we know there exists a constant \(C\), independent of \(h\), such that for all \((u,v) \in (V^2_h + H^1(\Omega)) \times V_h^2\),
    \[
    a^\sip_h(u,v) \leq C \| u \|_\sips \| v \|_\sip.
    \]
    Since \(u \in H^2(\Omega) + U^2_h \subset V^2_h + H^1(\Omega)\) and \(v \in V^0_h \subset V^2_h\), the inequality remains valid. 
\end{proof}

We are now ready to establish an a priori error estimates for Problem~\ref{eq:1d_dirichlet_disproblem}.

\begin{theorem}[Best approximation]\label{th:1d_dirichlet_approximation}
Let \(u\in H^2(\Omega)\) be the exact solution of~\eqref{eq:1d_dirichlet_poisson}, and let \(u_h\) be the discrete solution of~\eqref{eq:1d_dirichlet_disproblem}. Then there exists a constant \(C>0\), independent of \(h\), such that
\[
\|u-u_h\|_{\sip} \;\leq\; C \inf_{v_h\in U_h^2}\|u-v_h\|_{\sips}.
\]
\end{theorem}

\begin{proof}
Since \(u\) is continuous in \(\Omega\) and satisfies homogeneous Dirichlet boundary conditions, the jump terms vanish, i.e.,
\[
\jjump{u}|_F=0 \quad \forall F\in \mathcal{F}_h.
\]
Moreover, by definition of the discrete problem, 
\[
a_h^{\sip}(u-u_h,v)=0, \quad \forall v\in V_h^0.
\]

Using the inf--sup condition~\eqref{eq:1d_dirichlet_infsupcondition} together with boundedness~\eqref{eq:1d_dirichlet_boundedness}, we obtain
\[
C \|R^2u-u_h\|_{\sip}\,\|\Pi u-\Pi u_h\|_{\sip}
  \;\leq\; |a_h^{\sip}(R^2u-u_h,\Pi u-\Pi u_h)|
  \;=\; |a_h^{\sip}(R^2u-u,\Pi u-\Pi u_h)|.
\]
The last term is bounded by
\[
|a_h^{\sip}(R^2u-u,\Pi u-\Pi u_h)|
  \;\leq\; C \|u-R^2u\|_{\sips}\,\|\Pi u-\Pi u_h\|_{\sip}.
\]

Dividing both sides by \(\|\Pi u-\Pi u_h\|_{\sip}\) yields
\[
\|R^2u-u_h\|_{\sip} \;\leq\; C\|u-R^2u\|_{\sips}.
\]
Finally, the triangle inequality gives
\[
\|u-u_h\|_{\sip} 
 \;\leq\; \|u-R^2u\|_{\sip} + \|R^2u-u_h\|_{\sip} 
 \;\leq\; (1+C)\|u-R^2u\|_{\sips},
\]
which completes the proof.
\end{proof}

\begin{theorem}[Convergence rate]
Let \(u\in H^{3}(\Omega)\) be the exact solution of the model problem, and let \(u_h \in U_h^2\) be the numerical solution obtained from the proposed scheme. Then, for sufficiently small mesh size \(h\), there exists a constant \(C>0\), independent of \(h\), such that
\begin{align}
\|u-u_h\|_{\sips} &\leq C h^2, \label{1D_dirichlet_sip_error_estimate}\\
\|u-u_h\|_{L^2(\Omega)} &\leq C h^2. \label{1D_dirichlet_L2_error_estimate}
\end{align}
\end{theorem}

\begin{proof}
The estimate~\eqref{1D_dirichlet_sip_error_estimate} follows directly from the interpolation error bound in Theorem~\ref{th:interpolationestimate} combined with the a priori error estimate in Theorem~\ref{th:1d_dirichlet_approximation}.  

For the \(L^2\)-estimate, we recall from~\cite{daniele2012ipdg} the discrete inequality
\[
\|v\|_{L^2(\Omega)} \;\leq\; C \|v\|_{\sip}, \qquad \forall v \in V_h^2,
\]
where \(V_h^2\) denotes the standard DG space of degree \(2\). The estimate \eqref{1D_dirichlet_L2_error_estimate} follows immediately from~\eqref{1D_dirichlet_sip_error_estimate}.
\end{proof}
Although the $L^2$ error estimate in the above theorem appears suboptimal, numerical experiments in Sec.~\ref{sec:numericalresults} demonstrate that it is in fact sharp for this scheme, a behavior reminiscent of the second-order central scheme.

\subsubsection{2D Error Estimates}\label{subsec:2d_dirichlet_error}

The two-dimensional error analysis follows an approach analogous to the one-dimensional case. We first present several auxiliary lemmas that establish the properties of the reconstruction operator and equivalence of discrete norms, which are essential for proving stability and convergence.

\begin{lemma}\label{le:tensor}
	For any $R^k u \in \mathcal{Q}^{2m}(K)$, the following relations hold:
	\begin{equation}
		\begin{aligned}\label{eq:eq1}
			\int_{x_{j}}^{x_{j+1}}(R^ku)(x,y_j)\,\dx &= h\sum_{k=-m}^{m}\ccoef{k}{m}\bar{u}_{i+k+1/2,j+s+1/2},\\
			\int_{x_{j}}^{x_{j+1}}(R^ku)(x,y_{j+1})\,\dx &= h\sum_{k=-m}^{m}\ccoef{-k}{m}\bar{u}_{i+k+1/2,j+s+1/2},\\
			\int_{y_i}^{y_{i+1}}(R^ku)(x_i,y)\,\dy &= h\sum_{k=-m}^{m}\ccoef{k}{m}\bar{u}_{i+k+1/2,j+s+1/2},\\
			\int_{y_i}^{y_{i+1}}(R^ku)(x_{i+1},y)\,\dy &= h\sum_{k=-m}^{m}\ccoef{-k}{m}\bar{u}_{i+k+1/2,j+s+1/2}.
		\end{aligned}
	\end{equation}
	Similarly, for the derivatives, we have
	\begin{equation}
		\begin{aligned}\label{eq:eq2}
			\int_{x_{j}}^{x_{j+1}}(R^ku)_x(x,y_j)\,\dx &= \sum_{k=-m+1}^{m-1}\gcoef{k}{m} D_{(i+k,j+1/2)}u,\\
			\int_{y_i}^{y_{i+1}}(R^ku)_y(x_i,y)\,\dy &= \sum_{k=-m+1}^{m-1}\gcoef{k}{m} D_{(i+1/2,j+k)}u.
		\end{aligned}
	\end{equation}
\end{lemma}

The proof of Lemma~\ref{le:tensor} follows directly from the tensor-product structure of the reconstruction operator on the reference element $[0,1]^2$ and is analogous to the one-dimensional case. 
For brevity, the detailed derivation is omitted. To establish the equivalence between the finite volume method and the IPDG bilinear formulation in \eqref{eq:1d_dirichlet_disproblem}, we present the following lemma, which serves as the two-dimensional counterpart of Lemma~\ref{le:1d_grad_jump}.

\begin{lemma}\label{le:2d_grad_jump}
	In the two-dimensional case, for any $u \in V_h^0$, the reconstructed function $R^k u \in U_h^k$ satisfies
	\[
	\jjump{\nabla_h (R^k u)\!\cdot\! \mathbf{n}} = 0,
	\]
	for even $k$, where $\mathbf{n}$ denotes the unit outward normal vector on the cell interfaces.
\end{lemma}

We refer the reader to Appendix~\ref{ap:proof_of_le_2d_grad_jump} for the detailed proof.

We can also verify that the reconstructed solution satisfies the homogeneous Dirichlet boundary condition by prescribing appropriate cell averages on the virtual cells.

\begin{definition}
In two dimensions, the reconstruction operator $R^2$ is naturally defined on interior elements.  
For boundary elements $K$ with $\partial K \cap \partial \Omega \neq \varnothing$, the definition of the local reconstruction operator $R^2_K$ requires the use of cell averages on virtual elements.  

\medskip
\noindent\textbf{Edge elements.}  
Consider $\Gamma_{x=0}$ as an example. For $K=[0,h]\times [y_i,y_{i+1}],\; 1\leq i \leq N_2-1$, we seek $R^k_K u \in \mathcal{Q}^2(K)$ such that
\[
\int_{K'} R^k_K u \, \dx\,\mathrm{d}y \;=\; \int_{K'} u \, \dx\,\mathrm{d}y, 
\qquad \forall K'\subset S(K),
\]
with the virtual cell averages prescribed to enforce the boundary condition:
\[
R^k_K u(0,y_{i-1/2})=0, \quad R^k_K u(0,y_{i+1/2})=0, \quad R^k_K u(0,y_{i+3/2})=0.
\]

\medskip
\noindent\textbf{Corner elements.}  
For corner elements, additional constraints are imposed at the corner points.  
As an example, consider the bottom-left corner $K=[0,h]\times[0,h]$.  
We seek $R^k_K u \in \mathcal{Q}^2(K)$ such that
\[
\int_{K'} R^k_K u \, \dx\,\mathrm{d}y \;=\; \int_{K'} u \, \dx\,\mathrm{d}y, 
\qquad \forall K'\subset S(K), 
\]
with virtual cell averages prescribed to satisfy the boundary condition:
\[
R^k_K u(0,3h/2)=0, \; R^k_K u(0,h/2)=0, \; R^k_K u(0,0)=0, 
\; R^k_K u(h/2,0)=0, \; R^k_K u(3h/2,0)=0.
\]

\medskip
Analogous conditions are imposed for the other boundary edges and corners.
\end{definition}

Then, it can be easily show that the $R^2\bar{u}$ and consequently $u_h \in U_h^2$ satisfies the homogeneous Dirichlet boundary condition exactly, i.e. \[u_h|_{\partial \Omega}=0.\]
Similar to the 1D case, with the boundary condition strongly enforced in the reconstruction operator, we can establish the following theorem which establishes the equivalence of discrete norms.

\begin{lemma}[Norm equivalence]\label{2d_dirichlet_normequivalence}
	For the two-dimensional case, we define the semi-norm 
	\begin{align*}
		\vert u\vert_{1,h} = \left( \sum_{i=0}^{N}\sum_{j=0}^{N-1} (D_{(i,j+1/2)} u)^2 + \sum_{i=0}^{N-1}\sum_{j=0}^{N} (D_{(i+1/2,j)} u)^2 \right)^{1/2}.
	\end{align*}
	We show that $\vert u\vert_{1,h}$ is a norm in the $U^2_h$ space, denoted by $\Vert u\Vert_{1,h}$. Furthermore, there exist positive constants $C_1$ and $C_2$, such that
	\begin{equation}
		C_1 \Vert u\Vert_{\mathrm{sip}} \leq \Vert u\Vert_{1,h} \leq C_2 \Vert u\Vert_{\mathrm{sip}},
	\end{equation}
	where $C_1$ and $C_2$ are independent of $h$.
\end{lemma}

\begin{proof}
	To prove that $\vert\cdot \vert_{1,h}$ is a norm for $u \in U^2_h$, we first demonstrate its positive definiteness. The linearity and triangle inequality follow directly from standard arguments and will be omitted for brevity.

	Assuming 
	\begin{align*}
		\vert u\vert_{1,h} = 0,
	\end{align*}
	we have $\bar u \equiv \text{constant}$. From our treatment of the boundary condition, it follows that $\bar u = 0$, hence $u \equiv 0$. This concludes the positive definiteness of the semi-norm, confirming that it is indeed a norm, which we denote by $\Vert u \Vert_{1,h}$. 

	Using the previously established definitions, we compute $\Vert u \Vert_\sip$ and $\Vert u \Vert_{1,h}$ as follows:
	
	\medskip
	
	\noindent
For the sip norm, we have	
	\begin{equation*}
		\begin{aligned}
			\Vert u \Vert_{\sip}^2 
			&= \sum_{i=0}^{N-1} \sum_{j=0}^{N-1} 
			\left(
			\int_{x_i}^{x_{i+1}} \int_{y_j}^{y_{j+1}} (u_x^2 + u_y^2) \, \mathrm{d}x\,\mathrm{d}y\right)\\
		&\quad+\sum_{i=1}^{N-1}\sum_{j=0}^{N-1}\bigg(\frac{1}{h} \int_{y_j}^{y_{j+1}} \jjump{u(x_i,y)}^2 \, \mathrm{d}y\bigg)
			+\sum_{i=0}^{N-1}\sum_{j=1}^{N-1}\bigg(\frac{1}{h} \int_{x_i}^{x_{i+1}} \jjump{u(x,y_j)}^2 \, \mathrm{d}x
			\bigg) \\[4pt]
			&=\sum_{i,j} \left( \mathcal{R}^{(1)}_{i,j} + \mathcal{R}^{(2)}_{i,j} + \mathcal{R}^{(3)}_{i,j} + \mathcal{R}^{(4)}_{i,j} \right)\coloneq \mathcal{R}_1+ \mathcal{R}_2+ \mathcal{R}_3+ \mathcal{R}_4. \\[3pt]
		\end{aligned}
	\end{equation*}
\medskip
\noindent
The local terms are defined by
	\begin{equation*}
		\begin{aligned}
			\mathcal{R}^{(1)}_{i,j} &\coloneq \int_{x_i}^{x_{i+1}} \int_{y_j}^{y_{j+1}} u_x^2 \, \mathrm{d}x\,\mathrm{d}y, \\
			\mathcal{R}^{(2)}_{i,j} &\coloneq \int_{x_i}^{x_{i+1}} \int_{y_j}^{y_{j+1}} u_y^2 \, \mathrm{d}x\,\mathrm{d}y, \\
			\mathcal{R}^{(3)}_{i,j} &\coloneq \frac{1}{h} \int_{y_j}^{y_{j+1}} \jjump{u(x_i,y)}^2 \, \mathrm{d}y, \\
			\mathcal{R}^{(4)}_{i,j} &\coloneq \frac{1}{h} \int_{x_i}^{x_{i+1}} \jjump{u(x,y_j)}^2 \, \mathrm{d}x.
		\end{aligned}
	\end{equation*}

\medskip
\noindent
Likewise, for the discrete norm \(\Vert\cdot\Vert_{1,h}\), we obtain
		\begin{align*}
			\Vert u \Vert_{1,h}^2
		&= \sum_{i=0}^{N-1}\sum_{j=0}^{N} 
		\bigg(
		D_{(i+\frac{1}{2},j)} u
		\bigg)^2
		+ \sum_{i=0}^{N}\sum_{j=0}^{N-1} 
		\bigg(
		D_{(i,j+\frac{1}{2})} u
		\bigg)^2 \\[4pt]
		&\coloneq \sum_{i,j} \big( \mathcal{R}^{(5)}_{i,j} + \mathcal{R}^{(6)}_{i,j} \big)
		\coloneq \mathcal{R}_5+\mathcal{R}_6,
		\end{align*}
\noindent
here \begin{align*}
\mathcal{R}_{i,j}^{(5)}=\left(D_{{i+1/2,j}u}\right)^2,\\
\mathcal{R}_{i,j}^{(6)}=\left(D_{{i,j+1/2}u}\right)^2.
\end{align*}

To establish the equivalence of $\sipnorm{\cdot}$ and $\lnorm{\cdot}$, it suffices to show that $\mathcal{R}_1+\mathcal{R}_2$ and $\mathcal{R}_5+\mathcal{R}_6$ are mutually bounded, while $\mathcal{R}_3+\mathcal{R}_4$ is bounded by $\mathcal{R}_5+\mathcal{R}_6$ up to a constant $C$. 
	\begin{align*}
		&\text{Introduce local coefficients} \\
		&\quad
		\begin{alignedat}{2}
			A_{i,j} &\coloneq -\tfrac12 D_{(i,j-1/2)}u + \tfrac12 D_{(i+1,j-1/2)}u
			+ D_{(i,j+1/2)}u - D_{(i+1,j+1/2)}u \\[2pt]
			&\qquad\quad -\tfrac12 D_{(i,j+3/2)}u + \tfrac12 D_{(i+1,j+3/2)}u, \\
			B_{i,j} &\coloneq D_{(i,j-1/2)}u - 2D_{(i,j+1/2)}u + D_{(i,j+3/2)}u, \\
			C_{i,j} &\coloneq \tfrac12 D_{(i,j-1/2)}u - \tfrac12 D_{(i+1,j-1/2)}u
			-\tfrac12 D_{(i,j+1/2)}u + \tfrac12 D_{(i+1,j+1/2)}u, \\
			D_{i,j} &\coloneq -\tfrac16 D_{(i-1/2,j)}u + \tfrac16 D_{(i-1/2,j+1)}u
			-\tfrac{5}{12}D_{(i+1/2,j)}u + \tfrac{5}{12}D_{(i+1/2,j+1)}u \\[2pt]
			&\qquad\quad +\tfrac{1}{12}D_{(i+3/2,j)}u - \tfrac{1}{12}D_{(i+3/2,j+1)}u, \\
			E_{i,j} &\coloneq - D_{(i,j-1/2)}u + D_{(i,j+1/2)}u, \\
			F_{i,j} &\coloneq \tfrac13 D_{(i-1/2,j)}u + \tfrac56 D_{(i+1/2,j)}u - \tfrac16 D_{(i+3/2,j)}u.
		\end{alignedat}
		\\[4pt]
		&\text{Then}
		\qquad
		\mathcal{R}^{(1)}_{i,j}
		= \iint_{[0,1]^2} \big( A_{i,j}\,x y^2 + B_{i,j}\,x y + C_{i,j}\,y^2
		+ 2D_{i,j}\,x + E_{i,j}\,y + F_{i,j} \big)^2 \,\mathrm{d}x\,\mathrm{d}y.
	\end{align*}

	A direct computation together with Young's inequality yields:
	\begin{align*}
	\sum_{i=0}^{N-1}\sum_{j=0}^{N-1}\left(\mathcal{R}_{i,j}^{(1)}\right) \leq C\left(\sum_{i=0}^{N-1}\sum_{j=0}^{N-1}\mathcal{R}_{i,j}^{(5)}+\sum_{i=0}^{N-1}\sum_{j=0}^{N-1}\mathcal{R}_{i,j}^{(6)}\right).
	\end{align*}
	Similarly, we can derive that
		\begin{align*}
		\sum_{i=0}^{N-1}\sum_{j=0}^{N-1}\left(\mathcal{R}_{i,j}^{(2)}\right) \leq C\left(\sum_{i=0}^{N-1}\sum_{j=0}^{N-1}\mathcal{R}_{i,j}^{(5)}+\sum_{i=0}^{N-1}\sum_{j=0}^{N-1}\mathcal{R}_{i,j}^{(6)}\right),
	\end{align*}
	which implies 
	\begin{align*}
		\mathcal{R}_1+\mathcal{R}_2 \leq C  (\mathcal{R}_5+\mathcal{R}_6).
	\end{align*}

To establish the reverse inequality, we demonstrate that the quadratic forms associated with 
$\ddx{i}{j+{1}/{2}}u$ and $\ddx{i+{1}/{2}}{j}u$ are positive definite. 
Furthermore, the eigenvalues of the matrix corresponding to the quadratic form 
are uniformly positive constants independent of the indices $i,j$, the mesh size $N$, and the function $u$. Without loss of generality, we consider a representative element for illustration.

We next analyze the positivity of the quadratic form $\mathcal{R}_1+\mathcal{R}_2$. We define
\begin{align*}
	a_1 &= \ddx{i}{j-1/2}u, &\quad a_2 &= \ddx{i+1}{j-1/2}u, &
	a_3 &= \ddx{i}{j+1/2}u, &\quad a_4 &= \ddx{i+1}{j+1/2}u,\\
	a_5 &= \ddx{i}{j+3/2}u, &\quad a_6 &= \ddx{i+1}{j+3/2}u, &
	b_1 &= \ddx{i-1/2}{j}u, &\quad b_2 &= \ddx{i+1/2}{j}u,\\
	b_3 &= \ddx{i+3/2}{j}u, &\quad b_4 &= \ddx{i-1/2}{j+1}u, &
	b_5 &= \ddx{i+1/2}{j+1}u, &\quad b_6 &= \ddx{i+3/2}{j+1}u.
\end{align*}

	By direct observation, we find the following relationships:
\begin{align*}
	a_3 - a_1 &= b_3 - b_1, \quad
	a_6 - a_4 = b_6 - b_4, \\
	a_4 - a_2 &= b_5 - b_3, \quad
	a_5 - a_3 = b_4 - b_2,
\end{align*}
which allows us to eliminate $b_3, b_4, b_5, b_6$.
Substituting these relations into the definitions of $\mathcal{R}_{i,j}^{(1)}$ and $\mathcal{R}_{i,j}^{(2)}$, and simplifying, yields
	\[
\mathcal{R}_{i,j}^{(1)}+\mathcal{R}_{i,j}^{(2)} =  \bm{v}^\top \bm M \bm{v},
\qquad
\bm{v} = (a_1,a_2,a_3,a_4,a_5,a_6,b_1,b_2)^\top,
\]
	where \(\bm M\in\mathbb{R}^{8\times 8}\) is symmetric positive define whose eigenvalues are bounded above and below by positive constants independent of the mesh size. (see Appendix~\ref{appendix:spd_proof} for details). 

The positive definiteness of \(\bm M\) ensures the existence of a constant \(C\) such that 
\begin{align*}
	\mathcal{R}^{(1)}_{i,j}+\mathcal{R}^{(2)}_{i,j}=\bm v^\top\bm M\bm v\geq C\left(\sum_{i=1}^{6}a_i^2+b_1^2+b_2^2\right).
\end{align*}
Since $b_3, b_4, b_5, b_6$ are linear combinations of 
$a_1, \dots, a_6, b_1, b_2$, this estimate also controls all variables 
$a_i$ and $b_j$. Expressing the result in terms of the discrete derivatives
$D_{\bm{\theta}}u$, we obtain 
\begin{align*}
	\mathcal{R}_{5}+\mathcal{R}_{6} \le C(\mathcal{R}_{1}+\mathcal{R}_2).
\end{align*}

	
It remains to show that \(\mathcal{R}_3+\mathcal{R}_4\) can be bounded by \(\mathcal{R}_5+\mathcal{R}_6\)	up to a constant. We first prove that \(\mathcal{R}_3\) can be bounded by \(\Vert\cdot\Vert_{1,h}\), and the case for \(\mathcal{R}_4\) follows analogously. The analysis naturally separates into interior and boundary contributions.

	For the interior case, we have:
	\begin{align*}
		\mathcal{R}_{i,j}^{(3)} &= \int_{y_j}^{y_{j+1}} \left( -\frac{1}{12} \ddx{i}{j-3/2}u + \frac{1}{12} \ddx{i+1}{j-3/2}u + \frac{1}{4} \ddx{i}{j-1/2}u - \frac{1}{4} \ddx{i+1}{j-1/2}u \right. \\
		&\quad -\left. \frac{1}{4} \ddx{i}{j+1/2}u + \frac{1}{4} \ddx{i+1}{j+1/2}u + \frac{1}{12} \ddx{i}{j+3/2}u - \frac{1}{12} \ddx{i+1}{j+3/2}u \right)^2 \dy,
	\end{align*}
It follows that $\mathcal{R}^{(3)}_{i,j}$ is bounded by:
\[\mathcal{R}^{(3)}_{i,j} \leq C \left(\sum_{k=-m}^m\sum_{l=-m}^m (\ddx{i+k+1/2}{j+1}u)^2 +( \ddx{i+1}{j+k+1/2}u)^2\right),\]
here \(m=1\).

	For the boundary case, the jump term vanishes due to the boundary condition $u=0$ on $\partial \Omega$. Hence, a similar bound can be established.
Summing these local estimates over all indices $i$ and $j$, we obtain the global bound:
	\begin{align*}
		\mathcal{R}_{3} \leq C (\mathcal{R}_5+\mathcal{R}_6).
	\end{align*}
	Similarly, we have 
		\begin{align*}
		\mathcal{R}_{4} \leq C (\mathcal{R}_5+\mathcal{R}_6).
		\end{align*}


	
With both penalty terms bounded by $\Vert u \Vert_{1,h}^2$,  the equivalence of the $\Vert \cdot\Vert_{\sip}$ and $\Vert \cdot\Vert_{1,h}$ norms for $u \in U^2_h$ is a direct consequence.
	
\end{proof}

\begin{lemma}[Reconstruction boundedness] \label{lemma:stability_R2}
    Let $u \in V^0_h$ be such that its reconstruction $R^2u$ belongs to $U^2_h$. Then, the following stability estimate holds:
    \[
        \|u\|_{\sip} \leq C \|R^2 u\|_{\sip},
    \]
    where the constant $C > 0$ is independent of the mesh size $h$ and the function $u$.
\end{lemma}

\begin{proof}
    We begin by expanding the definition of the sip norm for $u \in V^0_h$:
    \begin{align*}
        \|u\|_{\sip}^2 
        &= \sum_{i=0}^{N}\sum_{j=0}^{N-1}\frac{1}{h}\int_{y_j}^{y_{j+1}} \jump{u}^2 \,dy+\sum_{i=0}^{N-1}\sum_{j=0}^{N}\frac{1}{h}\int_{x_i}^{x_{j+1}} \jump{u}^2 \,\dx.
    \end{align*}
    Since $u$ is a piecewise polynomial over the mesh, the integrals of the squared jumps on each edge can be evaluated exactly using a midpoint quadrature rule. This yields the discrete sum
    \begin{align*}
        \|u\|_{\sip}^2 
  &= \sum_{i=0}^{N}\sum_{j=0}^{N-1} \eta \, (\ddx{i}{j+1/2}u)^2 + \sum_{i=0}^{N-1}\sum_{j=0}^{N} \eta \, (\ddx{i+1/2}{j}u)^2.
    \end{align*}

    From Lemma~\ref{2d_dirichlet_normequivalence}, we have an explicit expression for $\|R^2 u\|_{\sip}^2$. This expression contains the same sum of squared jump terms as derived above, in addition to other non-negative terms originating from the element-wise gradients of the reconstructed function $R^2u$. Therefore, the norm of $u$ can be bounded by the norm of its reconstruction:
    \[
        \|u\|_{\sip}^2 \leq C \|R^2 u\|_{\sip}^2.
    \]
    Taking the square root of both sides concludes the proof.
\end{proof}

We now establish the inf--sup condition for the bilinear form \(a_h^{\sip}(\cdot,\cdot)\).

\begin{lemma}[Inf--sup condition]\label{the:2dinfsup}
	For the bilinear form \(a^\sip_h:U^2_h\times V^0_h\to \mathbb{R}\), the following inf--sup condition holds:
	\begin{align}
		\adjustlimits\inf_{u\in U_h^2\setminus\{0\}} 
		\sup_{v\in V_h^0\setminus\{0\}} 
		\frac{|a^\sip_h(u,v)|}{\Vert u\Vert_\sip \, \Vert v\Vert_\sip} \;\geq\; C,
		\label{eq:2d_dirichlet_infsup}
	\end{align}
	where \(C>0\) is a constant independent of $h$.
\end{lemma}

\begin{proof}
We prove the condition for the specific case $(r,k)=(0,2)$.	For \(u\in U^2_h,\, v\in V^0_h\), direct computation yields 
	\begin{align*}
a^\sip_h(u,v)&=\sum_{i=0}^{N-1}\sum_{j=1}^{N-1}\left(\sum_{k=-m+1}^{m-1}\ccoef{k}{m}\ddx{i+1/2}{j}u\right) \ddx{i+1/2}{j}v\\
&\quad+\sum_{i=1}^{N-1}\sum_{j=0}^{N-1}\left(\sum_{k=-m+1}^{m-1}\ccoef{k}{m}\ddx{i}{j+1/2}u\right)\ddx{i}{j+1/2}v\\
&\quad+\sum_{j=0}^{N}\left(\sum_{k=-m+1}^{m-1}\ccoef{k}{m}\ddx{0}{j+1/2}u\right)v_{1/2,j+1/2}-\sum_{j=0}^{N}\left(\sum_{k=-m+1}^{m-1}\ccoef{k}{m}\ddx{N}{j+1/2}u\right)v_{N-1/2,j+1/2}\\
&\quad+\sum_{i=0}^{N}\left(\sum_{k=-m+1}^{m-1}\ccoef{k}{m}\ddx{i+1/2}{0}u\right)v_{i+1/2,1/2}-\sum_{i=0}^{N}\left(\sum_{k=-m+1}^{m-1}\ccoef{k}{m}\ddx{i+1/2}{N}u\right)v_{i+1/2,N-1/2}\\
	&\quad+\eta\sum_{i=1,j=1}^{N,N-1}
\left(-\frac{1}{6}\ddx{i+1/2}{j-1}u+\frac{1}{3}\ddx{i+1/2}{j}u-\frac{1}{6}\ddx{i+1/2}{j+1}u\right)\ddx{i+1/2}{j}v \\
&\quad+\eta\sum_{i=1,j=1}^{N-1,N}
\left(-\frac{1}{6}\ddx{i-1}{j+1/2}u+\frac{1}{3}\ddx{i}{j+1/2}u-\frac{1}{6}\ddx{i+1}{j+1/2}u\right)\ddx{i}{j+1/2}v 
\end{align*}
		The enforced boundary condition gives:
	\begin{equation}\label{eq:boundary_ghost_points}
		\begin{aligned}
			&u_{j+1/2,-1/2} = -\frac{5}{2}u_{j+1/2,1/2} + \frac{1}{2}u_{j+1/2,3/2}, \quad &&j=0,\dots,N-1,\\
			&u_{j+1/2,N+1/2} = \frac{1}{2}u_{j+1/2,N-3/2} - \frac{5}{2}u_{j+1/2,N-1/2}, \quad &&j=0,\dots,N-1,\\
			&u_{-1/2,i+1/2} = -\frac{5}{2}u_{1/2,i+1/2} + \frac{1}{2}u_{3/2,i+1/2}, \quad &&i=0,\dots,N-1,\\
			&u_{N+1/2,i+1/2} = \frac{1}{2}u_{N-3/2,i+1/2} - \frac{5}{2}u_{N-1/2,i+1/2}, \quad &&i=0,\dots,N-1.
		\end{aligned}
	\end{equation}
We take \(v=\Pi u\). Substituting the equality \eqref{eq:boundary_ghost_points} into our expression yields: 
	\begin{align*}
	a^\sip_h(u,\Pi u)&=\sum_{i=0}^{N-1}\sum_{j=1}^{N-1}\left(\sum_{k=-m+1}^{m-1}\ccoef{k}{m}\ddx{i+1/2}{j}u\right) \ddx{i+1/2}{j}u\\
	&\quad+\sum_{i=1}^{N-1}\sum_{j=0}^{N-1}\left(\sum_{k=-m+1}^{m-1}\ccoef{k}{m}\ddx{i}{j+1/2}u\right)\ddx{i}{j+1/2}u\\
	&\quad+\sum_{j=0}^{N}\left(\sum_{k=-m+1}^{m-1}\ccoef{k}{m}\ddx{0}{j+1/2}u\right)\left(\tfrac{1}{3}\ddx{0}{j+1/2}u+\tfrac{1}{6}\ddx{1}{j+1/2}u\right)\\
	&\quad-\sum_{j=0}^{N}\left(\sum_{k=-m+1}^{m-1}\ccoef{k}{m}\ddx{N}{j+1/2}u\right)\left(\tfrac{1}{6}\ddx{N-1}{j+1/2}u+\tfrac{1}{3}\ddx{N}{j+1/2}u\right)\\
	&\quad+\sum_{i=0}^{N}\left(\sum_{k=-m+1}^{m-1}\ccoef{k}{m}\ddx{i+1/2}{0}u\right)\left(\tfrac{1}{3}\ddx{i+1/2}{0}u+\tfrac{1}{6}\ddx{i+1/2}{1}u\right)\\
	&\quad-\sum_{i=0}^{N}\left(\sum_{k=-m+1}^{m-1}\ccoef{k}{m}\ddx{i+1/2}{N}u\right)\left(\tfrac{1}{6}\ddx{i+1/2}{N-1}u+\tfrac{1}{3}\ddx{i+1/2}{N}u\right)\\
	&\quad+\eta\sum_{i=1,j=1}^{N,N-1}
	\left(-\tfrac{1}{6}\ddx{i+1/2}{j-1}u+\tfrac{1}{3}\ddx{i+1/2}{j}u-\tfrac{1}{6}\ddx{i+1/2}{j+1}u\right)\ddx{i+1/2}{j}u \\
	&\quad+\eta\sum_{i=1,j=1}^{N-1,N}
	\left(-\tfrac{1}{6}\ddx{i-1}{j+1/2}u+\tfrac{1}{3}\ddx{i}{j+1/2}u-\tfrac{1}{6}\ddx{i+1}{j+1/2}u\right)\ddx{i}{j+1/2}u.
\end{align*}

By Young's inequality and Lemma~\ref{lem:1d_deriv_representation}, it follows that
\begin{align*}
	a_h^{\sip}(u,\Pi u) \ge \frac{1}{h} \Bigg[
	& \sum_{i=0}^{N-1} \sum_{j=0}^{N} C_j \left( D_{i+1/2,j}u \right)^2 \\
	& + 
	\sum_{j=0}^{N-1} \sum_{i=0}^{N} C_i \left( D_{i,j+1/2}u \right)^2
	\Bigg]
\end{align*}
here
$$C_k =
\begin{cases}
	\dfrac{1}{3} - \dfrac{|1-\eta|}{12}, & \text{if } k \in \{0,N\},\\[6pt]
	\Big(1+\dfrac{\eta}{3}\Big) - \dfrac{|1-\eta|}{12} - \dfrac{|\eta|}{6}, & \text{if } k \in \{1,N-1\},\\[6pt]
	\Big(1+\dfrac{\eta}{3}\Big) - \dfrac{|\eta|}{3}, & \text{for } 2 \le k \le N-2.
\end{cases}$$
Consequently by choosing \(\eta\in (-\frac{3}{2},5)\), we derive that
\begin{align*}
a^{\sip}_h(u,\Pi u) \;&\geq\;C \Vert u\Vert_{1,h}\Vert \Pi u\Vert_{1,h}=C\Vert u\Vert_{1,h}\Vert R^2\Pi u \Vert_{1,h}
	\geq C\Vert u\Vert_\sip\Vert \Pi u\Vert_\sip.
\end{align*}
Therefore,
\[
\adjustlimits\inf_{u\in U^2_h\setminus\{0\}} \sup_{v\in V^0_h\setminus\{0\}} 
\frac{|a^\sip_h(u,v)|}{\Vert u\Vert_\sip \Vert v\Vert_\sip} \;\geq\; C,
\]
which completes the proof.
\end{proof}

The proof of boundedness for \(\bi\) is analogous to the one-dimensional case and is therefore omitted here.

\begin{lemma}[Boundedness]\label{le:2d_dirichlet_boundedness} 
The bilinear form \(a^\sip_h(\cdot,\cdot):(H^2+U^2_h)\times V^0_h\to \mathbb{R}\) is bounded.  
That is, there exists a constant \(C>0\), independent of \(h\), such that
\begin{equation}
	|a^\sip_h(u,v)| \;\leq\; C \Vert u\Vert_\sips \, \Vert v\Vert_\sip, 
	\qquad \forall u\in U^2_h,\ \forall v\in V^0_h.
	\label{eq:2d_dirichlet_boundedness}
\end{equation}
\end{lemma}

We are now ready to establish the a priori error estimates for the two-dimensional case of Problem~\eqref{eq:1d_dirichlet_disproblem}.

\begin{theorem}[Best approximation]
\label{thm:2d_energy_error}
Let \(u \in H^2(\Omega)\) be the exact solution of \eqref{eq:1d_dirichlet_poisson}, and let \(u_h \in U^2_h\) be the numerical solution of \eqref{eq:1d_dirichlet_disproblem}. Then there exists a constant \(C>0\), independent of \(h\), such that
\[
\Vert u - u_h \Vert_\sip \;\leq\; 
C \inf_{v_h \in U^2_h} \Vert u - v_h \Vert_\sips.
\]
\end{theorem}

\begin{proof}
Since \(u\) is continuous, the jump terms \(\jjump{u}\) vanish on all cell faces.  
Moreover, since both \(u\) and \(u_h\) satisfy the variational formulation, we have
\[
a^\sip_h(u - u_h, v) = 0, \qquad \forall v \in V^0_h.
\]

Applying the inf--sup condition~\eqref{eq:2d_dirichlet_infsup} together with the boundedness property~\eqref{eq:2d_dirichlet_boundedness}, we obtain
\begin{align*}
C \Vert R^2 u - u_h \Vert_\sip \, \Vert \Pi u - \Pi u_h \Vert_\sip
&\leq \big| a^\sip_h(R^2u - u_h, \Pi u - \Pi u_h) \big| \\
&= \big| a^\sip_h(R^2u - u, \Pi u - \Pi u_h) \big| \\
&\leq C' \, \Vert u - R^2u \Vert_\sips \, \Vert \Pi u - \Pi u_h \Vert_\sip.
\end{align*}
Dividing both sides by \(\Vert \Pi u - \Pi u_h \Vert_\sip\) gives
\[
\Vert R^2 u - u_h \Vert_\sip \;\leq\; C \Vert u - R^2u \Vert_\sips.
\]
Finally, by the triangle inequality,
\[
\Vert u - u_h \Vert_\sip 
\;\leq\; \Vert u - R^2 u \Vert_\sip + \Vert R^2 u - u_h \Vert_\sip
\;\leq\; (1+C)\, \Vert u - R^2 u \Vert_\sips.
\]
\end{proof}

\begin{theorem}[Convergence rate]\label{thm:2d_error_rates}
Let \(u \in H^{3}(\Omega)\) be the exact solution of the Poisson problem, and let \(u_h \in U^2_h\) be the numerical solution of the scheme.  
Then, for sufficiently small \(h\), there exists a constant \(C>0\) such that
\begin{align}
\Vert u - u_h \Vert_{\sips} &\leq C h^2, \label{2d_dirichlet_sip_error_estimate}\\
\Vert u - u_h \Vert_{L^2(\Omega)} &\leq C h^2. \label{2d_dirichlet_L2_error_estimate}
\end{align}
\end{theorem}

\begin{proof}
The energy-norm estimate \eqref{2d_dirichlet_sip_error_estimate} follows directly from Theorem~\ref{thm:2d_energy_error} and the interpolation properties of the reconstruction operator.  
For the \(L^2\)-estimate, we recall (see \cite{daniele2012ipdg}) that
\[
\Vert u \Vert_{L^2(\Omega)} \;\leq\; C \Vert u \Vert_\sip, 
\qquad \forall u \in V^2_h,
\]
where \(V^2_h\) denotes the discontinuous Galerkin space.  
Applying this inequality to the error function yields \eqref{2d_dirichlet_L2_error_estimate}.
\end{proof}

\subsection{Periodic Problem}\label{sec:periodic_ipdg}

In this section, we consider the periodic Poisson equation and briefly present the corresponding error estimates.
Since the periodic boundary condition can be readily enforced, the results for the periodic problem can be naturally extended to the case of even 
$k$. Most of the proofs follow the same arguments as those in the Dirichlet case~\eqref{eq:1d_dirichlet_poisson} and are therefore omitted for brevity.

We consider the Poisson problem: find \(u \in H^2(\Omega)\) such that
\begin{equation}\label{eq:poisson}
	-\Delta u = f, \quad  x \in \Omega,
\end{equation}
subject to periodic boundary conditions and the compatibility conditions
\[
\int_{\Omega} u \, \dx = 0, \; \int_{\Omega} f \, \dx = 0.
\]

To formulate the problem, we introduce a Lagrange multiplier to enforce the mean-zero condition.  
The corresponding bilinear form \(B(\cdot,\cdot)\) is defined for all \(u,v \in H^1(\Omega,\mathcal{T}_h)\) and \(\lambda,\mu \in \mathbb{R}\) by
\begin{align*}
B((u,\lambda),(v,\mu)) &= a^0_h(u,v) + b(\lambda,v) + b(u,\mu), 
\end{align*}
where
\begin{align*}
	a_h^{\sip}(u,v_h) &= -\sum_{K\in\mathcal{T}_h}\int_K \nabla_h u\!\cdot\!\nabla_h v_h \,\dx
	-\sum_{F\in\mathcal{F}_h}\int_F \aave{\nabla_h u}\!\cdot\! \mathbf{n}_F \,\jjump{v_h}\,\mathrm{d}s \\
	&\quad -\sum_{F\in\mathcal{F}_h}\int_F \jjump{u}\,\aave{\nabla_F v_h}\!\cdot\! \mathbf{n}_F \,\mathrm{d}s
\end{align*}
and
\[
b(u,\lambda) = \lambda \sum_{K\in\mathcal{T}_h} \int_K u \, \dx.
\]
Notice that we drop the penalty term for simplicity, since in the last section we found that the scheme \eqref{eq:1d_dirichlet_disproblem}  is well defined with $\eta = 0$. 

The discrete formulation is: find \((u_h,\lambda_h) \in U^k_h \times \mathbb{R}\) such that
\begin{equation}\label{eq:periodic_disproblem}
B((u_h,\lambda_h),(v_h,\mu_h)) = l(v_h),
\qquad \forall (v_h,\mu_h) \in  V^0_h \times \mathbb{R}.
\end{equation}
with $l(v) = \sum_{K\in\mathcal{T}_h}\int_K f v \, \dx.$

Let \(\bar{U}^k_h\) denote the subspace of \(U^k_h\) consisting of functions with zero mean value.
It can be shown that the discrete problem \eqref{eq:periodic_disproblem} is equivalent to the following reduced form:  
find \(u_h \in  \bar{U}^k_h\) such that
\[
a^0_h(u_h,v_h) = \tilde l(v_h), \qquad \forall v_h \in V^0_h,
\]
where
\[
\tilde l(v) = \sum_{K \in \mathcal{T}_h} \int_K (f - \bar f)\, v \, \dx, 
\qquad \tilde f = \frac{1}{|\Omega|}\int_\Omega f \, \dx.
\]

The periodic boundary condition can be naturally enforced within the virtual elements $\mathcal{T}_h^{\text{virt}}$ and extended to higher-order cases with larger virtual stencils. Therefore, the analysis for the periodic problem can be readily generalized to even values of $k$. 
The equivalence between the IPDG scheme~\eqref{eq:periodic_disproblem} and the finite volume method follows from Lemmas~\ref{le:1d_grad_jump} and~\ref{le:2d_grad_jump}, which implies that the second and third terms in~\eqref{eq:periodic_disproblem} vanish.

\subsubsection{1D Error Estimate}\label{subsec:1d_periodic_error}

In this subsection, we present the error estimates for the one-dimensional periodic problem.  Even though we will show the results for general high order,
since the proofs closely follow those of the Dirichlet case in Section~\ref{subsec:1d_dirichlet_error}, we only outline the key results to enhance clarity and readability.

\begin{lemma}[Norm equivalence]\label{eq:normequivalence}
Let \(\bar{U}^k_h\) be equipped with the semi-norm \(|\cdot|_{1,h}\) and the sip norm \(\|\cdot\|_{\sip}\) for even \(k\). 
Then \(|\cdot|_{1,h}\) defines a norm on \(\bar{U}^k_h\), denoted by \(\|\cdot\|_{1,h}\). 
Moreover, there exist constants \(C_1, C_2 > 0\), independent of \(h\), such that
\[
C_1\|u\|_{\sip} \le \|u\|_{1,h} \le C_2\|u\|_{\sip}, 
\qquad \forall\, u \in \bar{U}^k_h.
\]
\end{lemma}

\begin{lemma}[Reconstruction boundedness]\label{eq:operatorinjective}
	For any \(u \in V_h^0\), the reconstructed function \(R^k u \in U^k_h\) for even $k$ satisfies
	\[
	\|u\|_{\sip} \le C \|R^k u\|_{\sip},
	\]
	where \(C\) is independent of \(h\).
\end{lemma}

By Lemmas~\ref{eq:normequivalence} and~\ref{eq:operatorinjective}, the bilinear form \(a_h^0(\cdot,\cdot)\) satisfies the inf–sup and boundedness conditions.

\begin{lemma}[Inf--sup condition]\label{le:1d_periodic_infsupcondition}
	For the bilinear form \(a_h^0: U^k_h \times V^0_h \to \mathbb{R}\), there exists a constant \(C>0\), independent of \(h\) for even $k$, such that
	\[
	\adjustlimits\inf_{u\in U^k_h\setminus\{0\}} \sup_{v\in V^0_h\setminus\{0\}} 
	\frac{|a_h^0(u,v)|}{\|u\|_{\sip}\|v\|_{\sip}} \ge C.
	\]
\end{lemma}

\begin{lemma}[Boundedness]\label{le:1d_periodic_boundedness}
	The bilinear form \(a^0_h(\cdot,\cdot):(H^2+U^k_h)\times V^0_h\to\mathbb{R}\) is bounded for even $k$; that is, there exists a constant \(C>0\), independent of \(h\), such that
	\[
	|a^0_h(u,v)| \le C\|u\|_{\sips}\|v\|_{\sip}, 
	\qquad \forall u\in U^k_h,~\forall v\in V^0_h.
	\]
\end{lemma}

The proofs of Lemmas~\ref{le:1d_periodic_infsupcondition}-\ref{le:1d_periodic_boundedness} follow the same arguments as in the Dirichlet case and are thus omitted here for brevity.

\begin{theorem}[Best approximation]
	Let \(u \in H^{2}(\Omega)\) be the exact solution of the periodic Poisson problem~\eqref{eq:poisson}, and \(u_h \in U^k_h\) be the numerical solution of~\eqref{eq:periodic_disproblem} for even $k$. Then there exists a constant \(C>0\), independent of \(h\), such that
	\[
	\|u - u_h\|_{\sips} \le C \inf_{w_h \in U^k_h} \|u - w_h\|_{\sips}.
	\]
\end{theorem}

\begin{theorem}[Convergence rate]
	Let \(u \in H^{k+1}(\Omega)\) and \(u_h \in U^k_h\) be as above. Then, for sufficiently small \(h\),
	\begin{align}
		\|u - u_h\|_{\sips} &\le C h^k, \label{1D_periodic_sip_error_estimate}\\
		\|u - u_h\|_{L^2(\Omega)} &\le C h^k. \label{1D_periodic_L2_error_estimate}
	\end{align}
\end{theorem}

The proofs are analogous to those for the Dirichlet problem, following from the interpolation estimates in Theorem~\ref{th:interpolationestimate} and the stability properties established above, and are omitted for conciseness.

\subsubsection{2D Error Estimate}\label{subsec:2derror}

The two-dimensional error analysis follows the same framework as in the one-dimensional case. 
%
%
%
Based on Lemma~\ref{le:tensor}, we now establish the norm equivalence, which plays a central role in the subsequent error analysis.

\begin{lemma}[Norm equivalence]\label{2Dnormequivalence}
For the two-dimensional case,  for even $k$ define the semi-norm on $\bar{U}^k_h$ as
\[
|u|_{1,h}^2 = \sum_{i=0}^{N-1}\sum_{j=0}^{N-1} \left( D_{(i+1/2,j)} u \right)^2 + \sum_{i=0}^{N-1}\sum_{j=0}^{N-1} \left( D_{(i,j+1/2)} u \right)^2.
\]
Then $|\cdot|_{1,h}$ is a norm on $\bar{U}^k_h$, denoted by $\|\cdot\|_{1,h}$. 
Moreover, there exist constants $C_1, C_2 > 0$, independent of $h$, such that
\begin{equation}\label{eq:equi2}
C_1 \|u\|_{\sip} \le \|u\|_{1,h} \le C_2 \|u\|_{\sip}, \qquad \forall\, u \in \bar{U}^k_h.
\end{equation}
\end{lemma}
The proof follows the same argument as in the one-dimensional case, relying on the boundedness of the reconstruction operator and the consistency of the discrete gradients. 
Hence, it is omitted here for conciseness.

\begin{lemma}
For all $u \in V^0_h$ and the corresponding $R^k u \in U^k_h$ for even $k$, there exists a constant $C>0$, independent of $h$, such that
\[
\|u\|_{\sip} \le C \|R^k u\|_{\sip}.
\]
\end{lemma}

Next, we present the coercive and boundedness results of the bilinear form $a^0_h(\cdot, \cdot)$.
\begin{lemma}[Inf--sup condition]\label{the:2d_periodic_infsup}
For the bilinear form $a^0_h : U^k_h \times V^0_h \to \mathbb{R}$ for even $k$, the inf--sup condition holds:
\[
\adjustlimits\inf_{u\in U_h^k\setminus\{0\}} \sup_{v\in V_h^0\setminus\{0\}} 
\frac{|a^0_h(u,v)|}{\|u\|_{\sip}\|v\|_{\sip}} \ge C,
\]
where $C>0$ is a constant independent of $h$, $u$, and $v$.
\end{lemma}

\begin{lemma}[Boundedness]\label{le:2dboundedness}
The bilinear form $a^0_h(\cdot,\cdot): (H^{2}+U^k_h)\times V^0_h \to \mathbb{R}$ is bounded, i.e., there exists a constant $C>0$, independent of $h$, such that
\[
|a^0_h(u,v)| \le C \|u\|_{\sips}\|v\|_{\sip}, \qquad 
\forall\, u\in U^k_h,~ v\in V^0_h.
\]
\end{lemma}

\vspace{1em}

Finally, we present the a priori error estimates.

\begin{theorem}[Best approximation]
Let $u \in H^{2}(\Omega)$ be the exact solution of~\eqref{eq:poisson} and $u_h \in U^k_h$ the numerical solution of~\eqref{eq:periodic_disproblem}. 
Then there exists a constant $C>0$, independent of $h$, such that
\[
\|u - u_h\|_{\sip} \le C \inf_{v_h \in U^k_h} \|u - v_h\|_{\sips}.
\]
\end{theorem}

\begin{theorem}[Convergence rate]
Let $u \in H^{k+1}(\Omega)$ be the exact solution and $u_h \in U^k_h$ the numerical solution. 
Then, for sufficiently small $h$, there exists $C>0$ such that
\begin{align}
	\|u - u_h\|_{\sips} &\le Ch^k, \label{eq:2D_sip_error_estimate}\\
	\|u - u_h\|_{L^2(\Omega)} &\le Ch^k. \label{eq:2D_L2_error_estimate}
\end{align}
\end{theorem}

\section{Numerical Results}\label{sec:numericalresults}

In this section, we present a series of numerical experiments to validate the proposed scheme and examine its convergence behavior. 
The tests focus on both one- and two-dimensional Poisson problems with Dirichlet and periodic boundary conditions. 
Our primary goal is to verify that the observed convergence rates agree with the theoretical error estimates established in the previous sections.

\subsection{Convergence Study for the Dirichlet Problem}

\textbf{1D Test.}
We consider the one-dimensional Poisson equation
\begin{align*}
	-\Delta u &= f, \quad x \in [0,1],\\
	u(0) &= u(1) = 0.
\end{align*}
The exact solution is chosen as \(u(x) = x\sin(\pi x)\), from which the corresponding right-hand side \(f(x)\) is computed analytically.

\noindent \textbf{2D Test.}
We next consider the two-dimensional Poisson equation
\begin{align*}
	-\Delta u &= f, \quad (x,y) \in \Omega = [0,1]^2,\\
	u &= 0, \quad \text{on } \partial\Omega,
\end{align*}
with the exact solution \(u(x,y) = x\sin(\pi x)\,y\sin(\pi y)\). The source term \(f(x,y)\) is obtained accordingly.

\begin{figure}[H]
	\centering
	\includegraphics[width=0.42\textwidth]{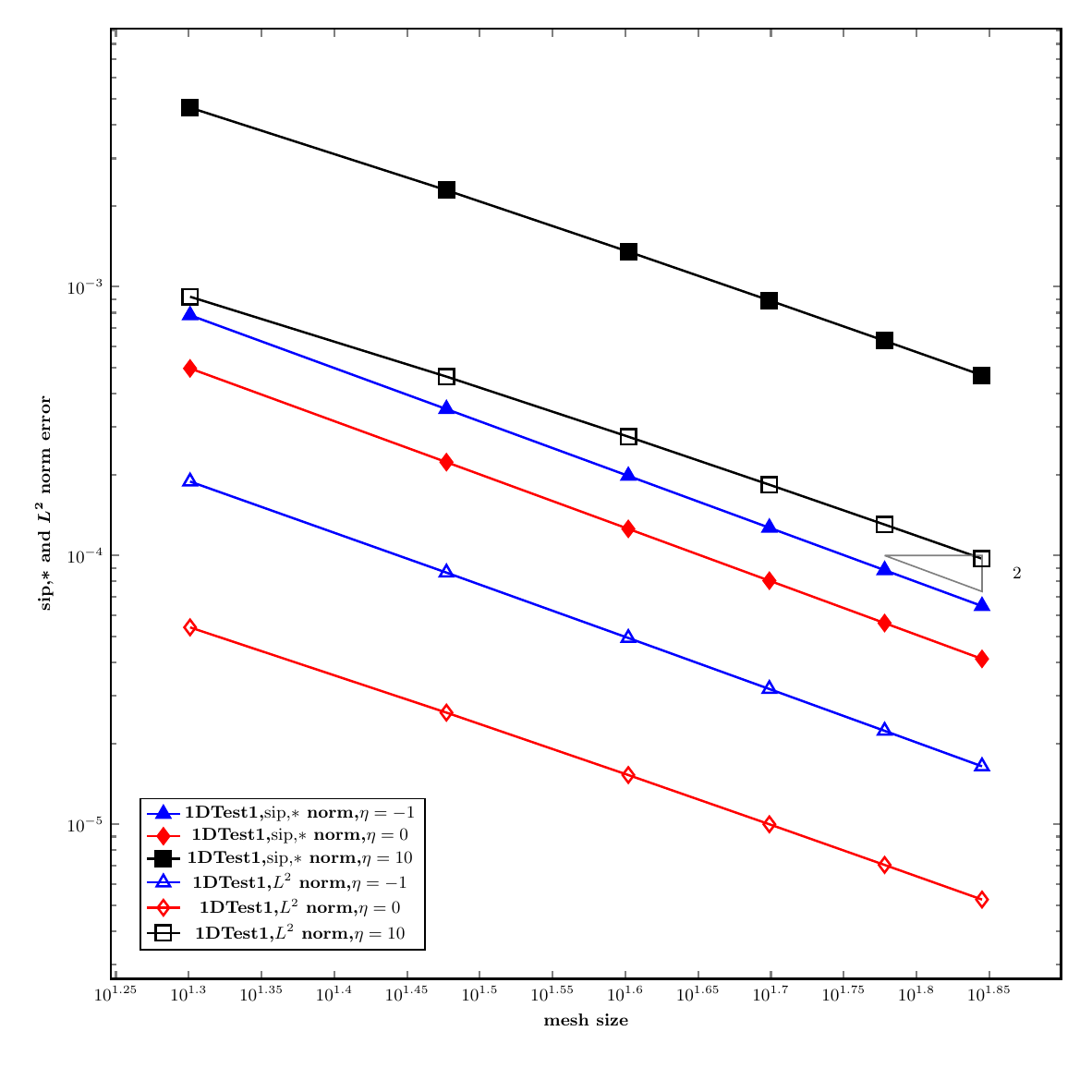}
	\hspace{0.08\textwidth}
	\includegraphics[width=0.42\textwidth]{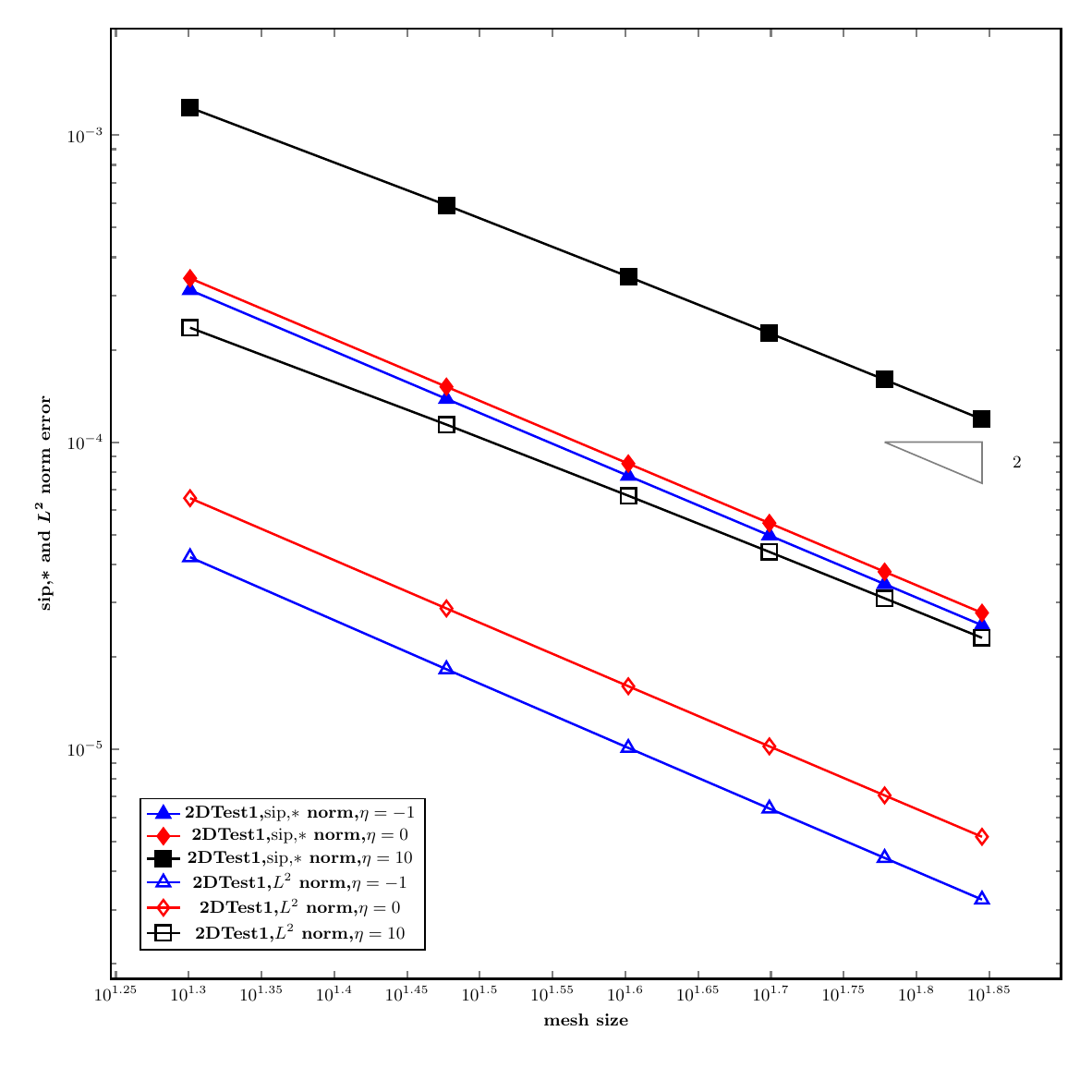}
	\caption{
	Convergence of the errors in the $\| \cdot \|_{L^2(\Omega)}$ and $\| \cdot \|_{\sips}$ norms for the Dirichlet problem in one dimension (left) and two dimensions (right).}
	\label{dirichlet_pic1}
\end{figure}

Figure~\ref{dirichlet_pic1} shows the convergence of the numerical errors measured in both $\| \cdot \|_{\sips}$ and $L^2$ norms. 
The results clearly demonstrate that the $\|u - u_h\|_{\sips}$ error converges at the optimal rate $\mathcal{O}(h^2)$, consistent with the theoretical prediction. 
Meanwhile, the $L^2$ error exhibits a suboptimal but consistent $\mathcal{O}(h^2)$ rate, confirming the stability and accuracy of the proposed finite volume formulation.

\subsection{Convergence Study for the Periodic Problem}

\textbf{1D Test.}
We consider the one-dimensional periodic Poisson equation
\begin{align*}
	-\Delta u &= f, \quad x \in [0,1],\\
	\int_0^1 u\,\dx &= 0, \qquad \int_0^1 f\,\dx = 0.
\end{align*}
The exact solution is taken as \(u(x) = \sin(2\pi x)\), and the corresponding \(f(x)\) is obtained by direct differentiation.

\noindent\textbf{2D Test.}
We further examine the two-dimensional periodic Poisson equation
\begin{align*}
	-\Delta u &= f, \quad (x,y)\in [0,1]^2,\\
	\int_{\Omega} u\,\dx\,\dy &= 0, \qquad \int_{\Omega} f\,\dx\,\dy = 0,
\end{align*}
with the exact solution \(u(x,y) = \sin(2\pi x)\sin(4\pi y)\). 
The corresponding right-hand side \(f(x,y)\) is again computed analytically.

\begin{figure}[H]
	\centering
	\includegraphics[width=0.42\textwidth]{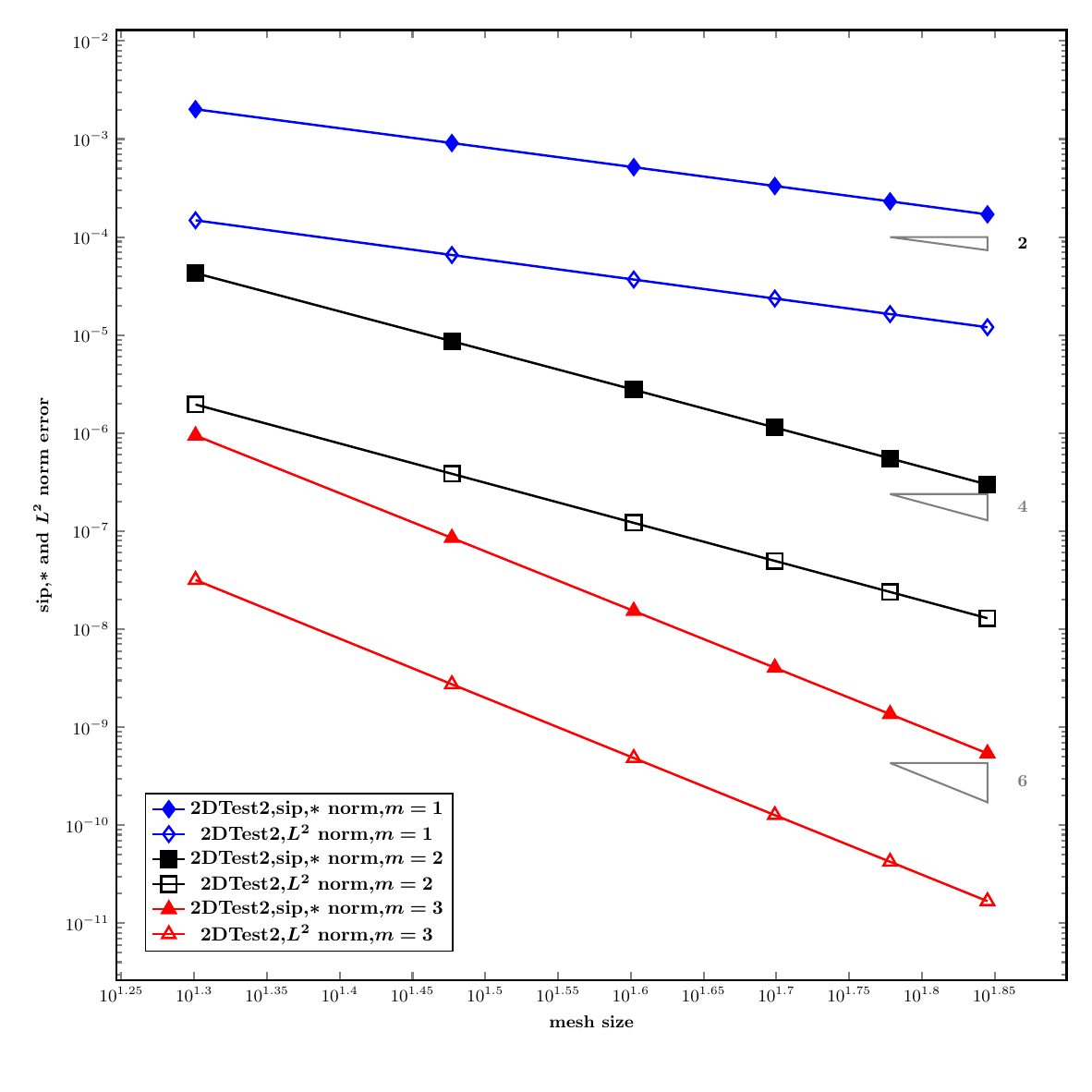}
	\hspace{0.08\textwidth}
	\includegraphics[width=0.42\textwidth]{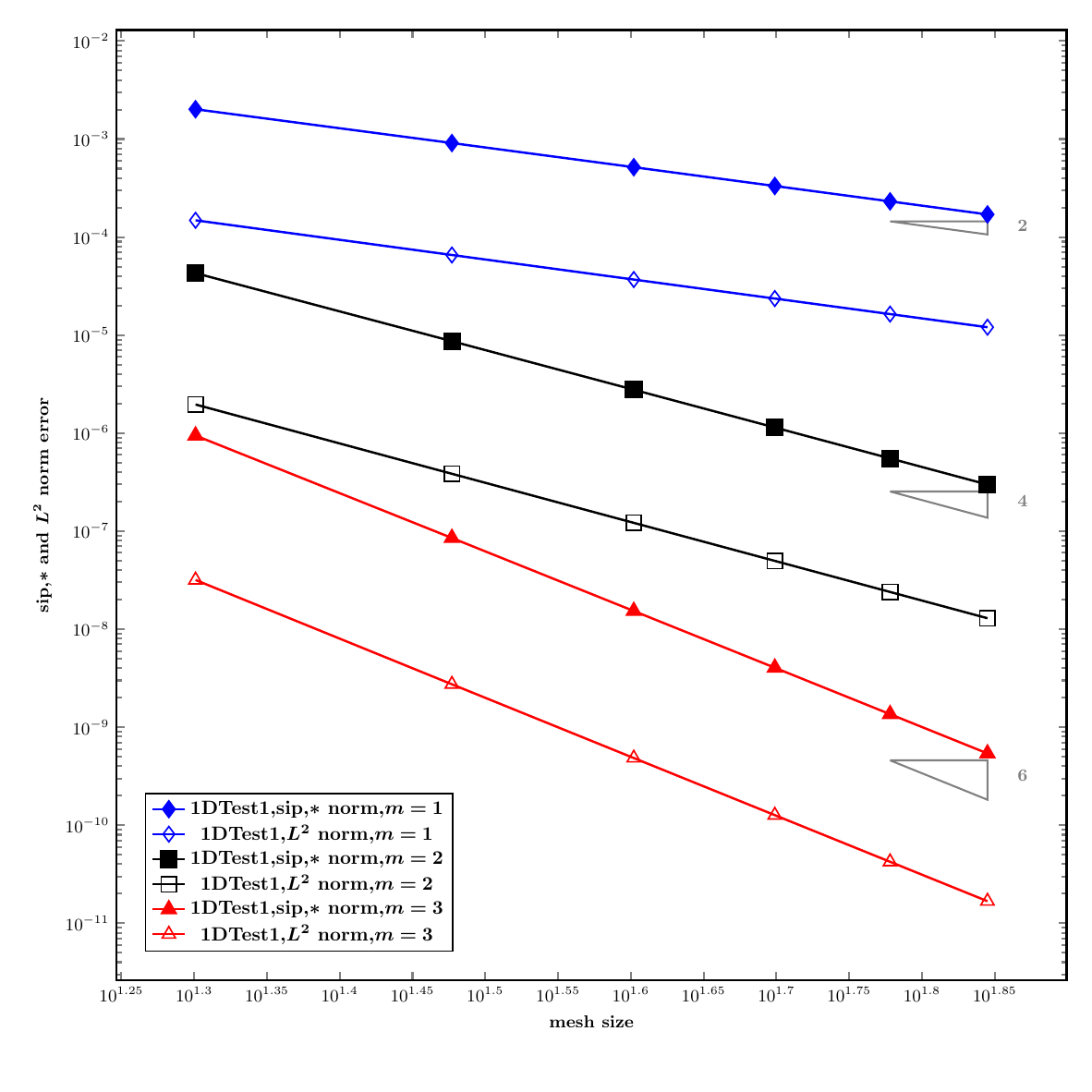}
	\caption{
	Convergence of the errors in the $\| \cdot \|_{L^2(\Omega)}$ and $\| \cdot \|_{\sips}$ norms for the periodic problem in one dimension (left) and two dimensions (right).}
	\label{periodic_pic}
\end{figure}

As shown in Figure~\ref{periodic_pic}, both the one- and two-dimensional periodic tests exhibit convergence rates in perfect agreement with the theoretical analysis. 
Specifically, the $\|u - u_h\|_{\sips}$ error converges at the optimal rate $\mathcal{O}(h^{k})$, while the $L^2$ error converges at the expected suboptimal rate $\mathcal{O}(h^{k})$. 
These numerical results confirm the validity of the theoretical error estimates and demonstrate the robustness of the proposed RDG-based finite volume formulation for both Dirichlet and periodic boundary conditions.

\section{Conclusions}\label{sec:conclusion}

In this work, we have developed and analyzed a high-order finite volume method for the Poisson problem based on the reduced discontinuous Galerkin (RDG) space.  
The key idea of the proposed approach is to adopt the piecewise constant space \(V_h^0\) as the test space and its corresponding higher-order RDG space \(U_h^k\) as the trial space, thereby formulating the scheme within a Petrov-Galerkin framework.  
This formulation naturally combines the local conservation property of finite volume methods with the high-order accuracy and flexibility of discontinuous Galerkin spaces, while maintaining a reduced number of degrees of freedom.

We have provided a rigorous theoretical analysis, establishing the stability and convergence of the method. In particular, we proved optimal-order convergence in the DG energy norm and suboptimal but consistent convergence in the \(L^2\) norm.  
Numerical experiments with both Dirichlet and periodic boundary conditions confirm the theoretical error estimates and demonstrate the robustness and efficiency of the proposed scheme.

The results presented here highlight the potential of the RDG-based Petrov-Galerkin framework as a bridge between classical finite volume and discontinuous Galerkin methods.  
This connection offers a promising foundation for the systematic development of high-order, locally conservative finite volume schemes with rigorous mathematical guarantees.  
Future work will extend this framework to nonlinear problems and time-dependent systems, where the advantages of RDG reconstruction and conservation properties can be further exploited.

\appendix

\section{Proof of Lemma~\ref{lem:1d_deriv_representation}}
\label{ap:proof_of_le:1d_deriv_representation}

The proof  utilizes the notation established during the proof of Lemma~\ref{le:1d_grad_jump} presented in Section~\ref{sec:dirichlet_ipdg}. 
For brevity, we refer the reader to that section for the definitions of \(\ccoef{s}{m}\), \( \left(L_k\right)_{k=-m}^{m+1} \). 
We begin by noting a key identity that follows directly from the property of the operator \( R^k \),
\begin{align}\label{eq:nature_identity}
\sum_{s=-m}^{m-1}\ccoef{s}{m}=0,
\end{align}
here 
\begin{align*}
\ccoef{s}{m}=\sum_{k=s}^{m}L''_{k+1}(0)=\sum_{k=s}^{m-1}L''_{k+1}(0),
\end{align*}
where the second equality holds since \(L_{m+1}(x)\) is an odd function.

In fact, recalling that 
\begin{align*}
	\widehat{(R^k u)_{\hat{x}}}\vert_{[0,1]}(0)=\sum_{s=-m}^{m-1}\ccoef{s}{m}\bu_{i+s+1/2},
\end{align*}
taking \(u\equiv 1\) yields~\eqref{eq:nature_identity}.

Then we have 
\begin{align}\label{eq:diff_expression}
\widehat{(R^k u)_{\hat{x}}}\vert_{[0,1]}(0)=\sum_{s=-m+1}^{m-1}\gcoef{s}{m}D_{(i+s)}u,
\end{align}
here
\begin{align*}
	\gcoef{s}{m}=\sum_{k=s}^{m-1}\ccoef{k}{m}.
\end{align*}

To prove the inequality~\eqref{eq:coeff_property}, we first 
establish four key structural properties of these \(\gcoef{k}{m}\) coefficients.

\begin{claim}[Symmetry]\label{claim:symmetry}
	The coefficients satisfy \(\gcoef{k}{m} = \gcoef{-k}{m}\) for \(k=1, \dots, m-1\).
\end{claim}
\begin{proof}[Proof of Claim~\ref{claim:symmetry}]
	To prove the symmetry property, it is sufficient to show that for $s=1,\dots,m-1$:
\[ \sum_{k=-s}^{s-1}\ccoef{k}{m}=0. \]
It is sufficient to show that for $s=1,\dots,m-1$:\[\ccoef{-s}{m}+\ccoef{s-1}{m}=0,\]
that is to prove that
\[\left(\sum_{k=-s}^{m-1}+\sum_{k=s-1}^{m-1}\right)L_{k+1}''(0)=0.\]

To demonstrate this, we define the polynomial
\[
S(x) := \sum_{k=-s}^{m-1}L_{k+1}(x)+\sum_{k=s-1}^{m-1}L_{k+1}(x).
\]
The desired identity is a direct consequence of the central symmetry of \(S(x)\) about the point \((0,\frac{1}{2})\). This symmetry implies that the second derivative, \(S''(x)\), is an odd function, and consequently, we must have \(S''(0) = 0\).
\end{proof}

\begin{claim}[Sign Alternation of \(L_k''(0)\)]\label{claim:sign1}
 \(L_{k}''(0)\) satisfies \(	\operatorname{sign}(L_{k}''(0))=(-1)^{k-1}\), for \(k=1, \dots, m-1\).
\end{claim}

\begin{proof}[Proof of Claim~\ref{claim:sign1}]
By definition, \(L_{k}(0)=0\). Hence, we can factor \(L_k(x)\) as
	\[L_k(x)=xg_k(x),\] here \(g_k(x)\) is a polynomial.
Differentiating twice, we obtain  \[L_k''(x)=2g_k'(x)+xg_k''(x),\]
and consequently,
	\begin{align}\label{eq:b1}
	L_k''(0)=2g_k'(0).
	\end{align}
To compute \(g_k(0)\), observe that  \(g_k(x)=\frac{L_k(x)}{x}\).
Taking logarithm gives
\begin{align}\label{eq:b2}\log{\vert g_k(x)\vert}=\log{\vert L_k(x)\vert}-\log{\vert x\vert}.\end{align}
Differentiating both sides yields
\begin{align}\label{eq:b3}
\frac{g_k'(x)}{g_k(x)}=\frac{L'_k(x)}{L_k(x)}-\frac{1}{x}=\sum_{\substack{r\in\{-m,\dots,m+1\}\\ r\neq 0,\, r\neq k}}\frac{1}{x-r}\coloneq H_k(x).
\end{align}
Therefore, the sign of \(L''_k(0)\) depends on the signs of \(g_k(0)\) and \(H_k(0)\).

For the explicit expression of \(L_k(x)\), we have 
	\[g_k(0)=\frac{\prod_{r=-m,\, r\neq 0,\,k}^{m+1}(-r)}{\prod_{r=-m,\, r\neq k}^{m+1}(k-r)},\]
it follows directly that 
\begin{align}\operatorname{sign}(g_k(0))=(-1)^{k-1}.\end{align}
	
	Next evaluating \(H_k(x)\) at \(x=0\) gives 
\begin{align}\label{eq:b4}H_k(0)=\sum_{\substack{r\in\{-m,\dots,m+1\}\\ r\neq 0,\, r\neq k}}\frac{1}{-r}=\frac{1}{k}-\frac{1}{m+1}>0\end{align}

Combining equations~\eqref{eq:b1},~\eqref{eq:b2},~\eqref{eq:b3} and~\eqref{eq:b4}, we conclude that \(L_k''(0)\) has the same sign as \((-1)^{k-1}\), which establishes the desired result.

\end{proof} 

\begin{claim}[Monotonicity]\label{claim:monotonicity}
	The function series satisfies $|L''_{k}(0)| > \lambda_k |L''_{k+1}(0)|$ for $k=1, \dots, m-2$, \(\forall \lambda_k < \frac{(k+1)^2(k+m+1)}{k^2(m-k)}\)
\end{claim}

\begin{proof}[Proof of Claim~\ref{claim:monotonicity}]
	Recalling the definition in the proof of~\ref{claim:sign1}, we have that
	\begin{align}\label{eq:exp}
		\vert L''_k(0)\vert=2\vert H_k(0)g_k(0)\vert&=2\left[\frac{m!(m+1)!}{k(k+m)!(m+1-k)!}\right]\cdot\left(\frac{1}{k}-\frac{1}{m+1}\right)\\
		&=\frac{2(m!)^2}{k^2(m+k)!(m-k)!}.
	\end{align}
	
To verify monotonicity, consider the ratio
\begin{align*}
	\frac{\vert L''_{k+1}(0)\vert}{\vert L''_{k}(0)\vert}&= \frac{k^2 (k+m)! (m-k)!}{(k+1)^2 (k+m+1)! (m-k-1)!}\\
	&=\frac{k^2(m-k)}{(k+1)^2(k+m+1)}<\frac{1}{\lambda_k},
\end{align*}	
equivalently,
\begin{align*}
	|L''_{k}(0)| > \lambda_k |L''_{k+1}(0)| .
\end{align*}
%

\end{proof}

\begin{claim}[Sign Alternation of \(\gcoef{k}{m}\)]\label{claim:sign2}
 	The coefficients satisfy \(\operatorname{sign}(\gamma_{k,m}) = (-1)^k\) for \(k=1, \dots, m-1\).
\end{claim}

\begin{proof}[Proof of Claim~\ref{claim:sign2}]
	
	We first analyze the coefficients $\ccoef{s}{m}$. The sum $\ccoef{s}{m}$ is an alternating series, as established by Claim~\ref{claim:sign1}. Furthermore, Claim~\ref{claim:monotonicity} implies that the terms of this series are strictly decreasing in magnitude (i.e., $|L_k''(0)| > |L_{k+1}''(0)|$, as previously shown). By the properties of such an alternating series, its sign is determined by the sign of its first term. Therefore,$$\operatorname{sign}\big(\ccoef{s}{m}\big) = \operatorname{sign}\big(L_{s+1}''(0)\big) = (-1)^{(s+1)-1} = (-1)^{s}.$$
	We now prove that for \(s=0,\dots, m-1\),
	\begin{align*}
		\vert \ccoef{s}{m}\vert>\vert \ccoef{s+1}{m}\vert.
	\end{align*}
We distinguish between the cases where \(s\) is even and where \(s\) is odd.
	
(i) For the even case \(s=2s'\), the goal is to establish
\begin{align}\label{eq:sum1}
	\ccoef{2s}{m}+\ccoef{2s+1}{m}>0.
\end{align}
To this end, the analysis is further divided according to the parity of \(m\).

(a) When \(m\) is odd, the argument proceeds by backward induction on \(s'\).  
 
For the base case \(s'=\frac{m-1}{2}\),
\begin{align*}
\ccoef{m-1}{m}+\ccoef{m}{m}=L_{m}''(0)>0.
\end{align*}
Suppose that for some \(s'+1\leq \frac{m+1}{2}\), the following inequality holds  
\begin{align*}
\ccoef{2s'+2}{m}+\ccoef{2s'+3}{m}>0.	
\end{align*}
It remains to verify that it also holds for \(s'\).  
Using the recursive relation, one obtains
\begin{align*}
&(\ccoef{2s'+2}{m}+\ccoef{2s'+3}{m})-(\ccoef{2s'}{m}+\ccoef{2s'+1}{m})\\
&=2L_{2s'+3}''(0)+L''_{2s'+2}(0)+L''_{2s'+4}(0).
\end{align*}
Applying Claim~\ref{claim:sign1} yields
\begin{align*}
&\quad 2L_{2s'+3}''(0)+L''_{2s'+2}(0)+L''_{2s'+4}(0)\\
&=-\vert L''_{2s'+2}(0)\vert+\vert 2L_{2s'+3}''(0)\vert-\vert L''_{2s'+4}(0)\vert.
\end{align*}
By~\eqref{eq:exp} in Claim~\ref{claim:monotonicity}, a straightforward calculation shows that
\begin{align*}
	(\ccoef{2s'+2}{m}+\ccoef{2s'+3}{m})-(\ccoef{2s'}{m}+\ccoef{2s'+1}{m})<0,
\end{align*}
which preserves the sign, completing the backward induction.

\smallskip

(b) When \(m\) is even, the proof follows the same inductive strategy.   For the base case \(s'=\frac{m}{2}-1\),
\begin{align*}
	\ccoef{m-2}{m}+\ccoef{m-1}{m}=L_{m-1}''(0)+2L_{m}''(0).
\end{align*}
Since \(\lambda_{m-1}=\tfrac{2m^3}{(m-1)^2}>2\), it follows that
\begin{align*}
	\ccoef{m-2}{m}+\ccoef{m-1}{m}=L_{m-1}''(0)+2L_{m}''(0)>0.
\end{align*}

The inductive step proceeds in exactly the same manner as in the odd-\(m\) case, thereby completing the proof for even \(m\).

(ii) For the odd case \(s=2s'-1\), an entirely analogous argument shows that 
\begin{align}\label{eq:sum2}
	\ccoef{2s'-1}{m}+\ccoef{2s'}{m}<0.
\end{align}
Combining~\eqref{eq:sum1} and~\eqref{eq:sum2}, we conclude that 
\(\lvert \ccoef{s}{m}\rvert\) is strictly decreasing in \(s\).

\smallskip

The same reasoning applies to the sign property of \(\gcoef{s}{m}\), yielding the identical conclusion.  
This completes the proof.

\end{proof}

	
With properties~\ref{claim:symmetry},~\ref{claim:sign1},~\ref{claim:monotonicity}, and~\ref{claim:sign2} established above, 
we are now ready to prove~\eqref{eq:coeff_property}.  Recalling the definition of \(\gcoef{k}{m}\) and \(\ccoef{k}{m}\), 
\begin{align*}
\gcoef{0}{m}=\sum_{k=0}^{m-1}\ccoef{k}{m}&=\sum_{k=0}^{m-1}\sum_{s=k}^{m-1}L_{s+1}''(0)\\
&=\sum_{s=0}^{m-1}(s+1)L_{s+1}''(0).
\end{align*}

Furthermore,
\begin{align*}
	\sum_{k=-m,\, k\neq 0}^{m-1}\vert\gcoef{k}{m}\vert&=2\sum_{k=1}^{m-1}(-1)^k\sum_{s=k}^{m-1}\ccoef{s}{m}\\
	&=2\sum_{k=1}^{m-1}\sum_{s=k}^{m-1}\sum_{r=s}^{m-1}(-1)^kL_{r+1}''(0)\\
	&=-2 \sum_{r=1}^{m-1} \left\lceil \frac{r}{2} \right\rceil L_{r+1}''(0).
\end{align*}

Hence,
\begin{align}\label{eq:coeff_property_minus}
\gcoef{0}{m}-\sum_{k=-m,\, k\neq 0}^{m-1}\vert\gcoef{k}{m}\vert = \sum_{j=1}^{m} C_j L_j''(0),
\end{align}
where the coefficients $C_j$ are defined by
\[
C_j =
\begin{cases}
	2j - 1, & \text{if $j$ is odd},\\[3pt]
	2j, & \text{if $j$ is even.}
\end{cases}
\]

Expanding the right-hand side of~\eqref{eq:coeff_property_minus} and Using Claim~\ref{claim:monotonicity} gives the following two cases.
When $m$ is even, both even and odd indices contribute positively, leading to
\begin{align*}
	\gcoef{0}{m}-\sum_{k=-m,\, k\neq 0}^{m-1}\vert\gcoef{k}{m}\vert=&\sum_{l=1}^{m/2}C_{2l}L''_{2l}(0)+\sum_{l=1}^{m/2}C_{2l-1}L''_{2l-1}(0)\\
	&=\sum_{l=1}^{m/2}4lL''_{2l}(0)+\sum_{l=1}^{m/2}(4l-3)L''_{2l-1}(0)> 0.
\end{align*}

When $m$ is odd, the same argument yields

\begin{align*}
	\gcoef{0}{m}-\sum_{k=-m,\, k\neq 0}^{m-1}\vert\gcoef{k}{m}\vert=&\sum_{l=1}^{(m+1)/2}C_{2l}L''_{2l}(0)+\sum_{l=1}^{(m+3)/2}C_{2l-1}L''_{2l-1}(0)\\
	&=\sum_{l=1}^{(m-1)/2}4lL''_{2l}(0)+\sum_{l=1}^{(m+1)/2}(4l-3)L''_{2l-1}(0)\\
	&>\sum_{l=1}^{(m-1)/2}4lL''_{2l}(0)+\sum_{l=1}^{(m-1)/2}(4l-3)L''_{2l-1}(0)\geq 0.
\end{align*}

Hence, the desired property~\eqref{eq:coeff_property} holds for all $m \ge 2$.

\section{Proof of  Lemma~\ref{le:2d_grad_jump}}
\label{ap:proof_of_le_2d_grad_jump}
\begin{proof}[Proof of Lemma~\ref{le:2d_grad_jump}]
	We prove the result for the boundary face $\{x_i\} \times [y_j, y_{j+1}]$; the proofs for other faces are analogous and thus omitted for brevity.  Let \(\{\phi_p(x)\}_{p=-m}^m\) and \(\{\phi_q(y)\}_{q=-m}^m\) denote the set of one-dimensional basis functions on reduced Galerkin space. 
	For each element \(K_{i,j}=[x_j,x_{j+1}]\times[y_i,y_{i+1}]\), there exists an affine mapping
	from the reference element \(\hat{K}=[0,1]^2\). This mapping, \(F_{i,j}:\hat K\to K_{i,j}\), transforms the reference coordinates \(\hat x=(\hat x,\hat y)\in \hat K\) to the physical coordinate 
	\(x=(x,y)\in K_{i,j}\) as follow:
	\begin{align*}
		{x} = F_{ij}(\hat{{x}}) =
		\begin{pmatrix}
			x_j + \hat{x} (x_{j+1} - x_j) \\
			y_i + \hat{y} (y_{i+1} - y_i)
		\end{pmatrix}
		=
		\begin{pmatrix}
			x_j + \hat{x} h_{j} \\
			y_i + \hat{y} h_{i}
		\end{pmatrix},
	\end{align*}
	where $h_j = x_{j+1} - x_j$ and $h_i = y_{i+1} - y_i$ are the element widths in the $x$ and $y$ directions, respectively.
		The derivative of the reconstructed function $\widehat{R^k u}$ from Lemma~\ref{le:tensor} is
	\begin{align}\label{eq:reconstruction_expression}
	(\widehat{R^k u})_{\hat x}(\hat x, \hat y)
	= \sum_{p=0}^{2m}\sum_{q=0}^{2m} c_{pq} \, (\widehat\phi_p)'(\hat x)\, \widehat\phi_q(\hat y).
	\end{align}
	
	By the property of the reconstruction operator \(R^k\), for $ s, k = -m, \dots, m-1$,
	\begin{align}\label{eq:reconstruction_integration}
		\bu_{i+k+1/2,j+s+1/2}=\int_{k}^{k+1}\int_{s}^{s+1}\hat{R^k u}(\hat x,\hat y)\,\mathrm{d}\hat x\mathrm{d}\hat y,
	\end{align}
substituting~\eqref{eq:reconstruction_expression} and~\eqref{eq:reconstruction_integration}, we derive that
	\begin{align*}
	\bu_{i+k+1/2,j+s+1/2}=\sum_{p=0}^{2m}\sum_{q=0}^{2m}c_{pq}\int_{k}^{k+1}\hat\phi_p(\hat x)\,\mathrm{d}\hat x\int_{s}^{s+1}\hat{\phi}(\hat y)\,\mathrm{d}\hat y.
\end{align*}

Define $\bm{A} = (a_{p,k})$ where $0 \le p \le 2m$ and $-m \le k \le m$, here \(a_{p,k}\) is defined as follow:
\[
a_{p,k}=\int_{k}^{k+1}\hat{\phi}_p(\hat x)\,\mathrm{d}\hat x,
\]
and \(k\) denotes the cell index and \(p\) the polynomial degree. By definition, \(\bm A\) is invertible.

Let \(\bm U\) denote \((\bu_{i+k+1/2,j+s+1/2})\) where $-m \le k \le m$ and $-m \le s \le m$, and \(\bm C\) denote  \((c_{pq})\)   where $0 \le p \le 2m$ and $0 \le q \le 2m$. Then we have \(\bm C=(\bm A^\top)^{-1}\bm U(\bm A)^{-1}\).

Evaluating at $\hat x = 0$ gives
	\begin{align*}
		(\widehat{R^k u})_{\hat x}(0, \hat y)
		&= \sum_{p=0}^{2m}\sum_{q=0}^{2m} c_{pq} (\widehat\phi_p)'(0) \widehat\phi_q(\hat y)
		= ((\widehat\phi)'_p(0))^\top \bm C (\widehat\phi_q(\hat y)) \\
		&= ((\widehat\phi)'_p(0))^\top (\bm A^\top)^{-1} \bm U \bm A^{-1} (\widehat\phi_q(\hat y)).
	\end{align*}
	
	Our goal is to show that this expression is equivalent to
	\begin{align}\label{eq:grad}
		(\widehat{R^k u})_{\hat x}(0, \hat y)
		= \sum_{q=0}^{2m} \Bigl[\sum_{k=-m+1}^{m-1} c_m^k D_{(i+q+1/2, j+k)}u\Bigr]^\top (\bm A)^{-1} \widehat\phi_q(\hat y).
	\end{align}
	
	To prove~\eqref{eq:grad}, it suffices to show that for $q = 0, \dots, 2m$,
	\[
	((\bm A^\top)^{-1}\bm u_{j+q+1/2})^\top (\widehat\phi'_p(0))
	= \sum_{k=-m+1}^{m-1} c_m^k D_{(i+q+1/2, j+k)}u,
	\]
	here \(\bm u_{j+q+1/2}\) denotes the column vector of \(\bm U\) corresponding to the \text{y-index} \(j+q+1/2\).
	In the one-dimensional case, it holds that
	\[
	((\bm A^\top)^{-1}\tilde{\bm u})^\top (\widehat\phi'_p(0))
	= \sum_{k=-m+1}^{m-1} c_m^k D_{(i+k)}u,
	\]
	where $\tilde{\bm u}$ is the column vector of one-dimensional cell averages.
	This establishes~\eqref{eq:grad}.
	
	Consequently,
	\[
	((\widehat{R^k u})_{\hat x}|_{[0,1]^2})(0, \hat y)
	= \sum_{q=0}^{2m} \Bigl[\sum_{k=-m+1}^{m-1} c_m^k D_{(i+q+1/2, j+k)}u\Bigr]^\top (\bm A)^{-1} \widehat\phi_q(\hat y).
	\]
	Similarly, for the adjacent element $[-1,0]\times[0,1]$, we obtain
	\[
	((\widehat{R^k u})_{\hat x}|_{[-1,0]\times[0,1]})(0, \hat y)
	= \sum_{q=0}^{2m} \Bigl[\sum_{k=-m+1}^{m-1} c_m^k D_{(i+q+1/2, j+k)}u\Bigr]^\top (\bm A)^{-1} \widehat\phi_q(\hat y).
	\]
	Since the coefficients $\ccoef{k}{m}$ are uniquely determined by Lemmas~\ref{le:1d_grad_jump} and~\ref{lem:1d_deriv_representation}, the two expressions coincide, leading to
	\[
	\jjump{(R^k u)_x} = 0 \quad \text{on } \{x_i\} \times [y_j, y_{j+1}],
	\]
	which completes the proof.
\end{proof}

\section{Proof of the symmetric positive definiteness of the local matrix}\label{appendix:spd_proof}

In this appendix, we verify that the local matrix $M$ introduced in Lemma~\ref{2d_dirichlet_normequivalence} 
is symmetric positive definite. All notations and indexing conventions follow those in the proof of Lemma~\ref{2d_dirichlet_normequivalence}.

By direct computation, the matrix $M$ is given by
\begin{align*}
	M = \frac{1}{1080} \left[
	\begin{array}{cccc cc cc}
		47/135 & 11/1080 & -37/216 & -1/270 & -191/1080 & -7/1080 & -1/3 & -1/6 \\
		11/1080 & 2/135 & -1/270 & -1/216 & -7/1080 & -11/1080 & 0 & 0 \\
		-37/216 & -1/270 & 73/108 & 47/270 & -37/216 & -1/270 & 1/6 & -1/6 \\
		-1/270 & -1/216 & 47/270 & 37/108 & -1/270 & -1/216 & 0 & 0 \\
		-191/1080 & -7/1080 & -37/216 & -1/270 & 47/135 & 11/1080 & 1/6 & 1/3 \\
		-7/1080 & -11/1080 & -1/270 & -1/216 & 11/1080 & 2/135 & 0 & 0 \\
		-1/3 & 0 & 1/6 & 0 & 1/6 & 0 & 1/3 & 1/6 \\
		-1/6 & 0 & -1/6 & 0 & 1/3 & 0 & 1/6 & 1/3
	\end{array}
	\right].
\end{align*}

All eigenvalues of $M$ are strictly positive and independent of the mesh parameter $N$, namely,
\begin{align*}
	[\lambda_i]_{i=1}^8 = \frac{1}{1080}
	\begin{bmatrix}
		26/21477 \\[2pt]
		52/4825  \\[2pt]
		66/5381  \\[2pt]
		202/3175 \\[2pt]
		209/1131 \\[2pt]
		435/643  \\[2pt]
		1550/839 \\[2pt]
		1266/625
	\end{bmatrix}.
\end{align*}
Therefore, $M$ is symmetric positive definite with eigenvalues uniformly bounded away from zero and infinity, independent of $N$, which completes the proof.

 \section*{Declarations}

 \subsection*{Data Availability}
 The datasets generated during the current study are available from the corresponding author upon reasonable request. They support
 our published claims and comply with field standards.
 
  \subsection*{Competing of interest}
The authors declare that they have no known competing financial interests or personal relationships that could have appeared to influence the work reported in this paper.

\bibliographystyle{abbrv}
\bibliography{ref_4}

@article {jones2000fvm,
	AUTHOR = {Jones, W. P. and Menzies, K. R.},
	TITLE = {Analysis of the cell-centred finite volume method for the
	diffusion equation},
	JOURNAL = {J. Compute. Phys.},
	FJOURNAL = {Journal of Computational Physics},
	VOLUME = {165},
	YEAR = {2000},
	NUMBER = {1},
	PAGES = {45--68},
	ISSN = {0021-9991,1090-2716},
	MRCLASS = {65M70 (65M12)},
	MRNUMBER = {1795392},
	DOI = {10.1006/jcph.2000.6595},
	URL = {https://doi.org/10.1006/jcph.2000.6595},
}

@article {liu1994weno,
	AUTHOR = {Liu, Xu-Dong and Osher, Stanley and Chan, Tony},
	TITLE = {Weighted essentially non-oscillatory schemes},
	JOURNAL = {J. Compute. Phys.},
	FJOURNAL = {Journal of Computational Physics},
	VOLUME = {115},
	YEAR = {1994},
	NUMBER = {1},
	PAGES = {200--212},
	ISSN = {0021-9991,1090-2716},
	MRCLASS = {65M99},
	MRNUMBER = {1300340},
	DOI = {10.1006/jcph.1994.1187},
	URL = {https://doi.org/10.1006/jcph.1994.1187},
}

@article{harten1987uniformly,
	title = {Uniformly high order accurate essentially non-oscillatory schemes, {III}},
	journal = {J. Comput. Phys.},
	fjournal = {Journal of Computational Physics},
	volume = {71},
	number = {2},
	pages = {231-303},
	year = {1987},
	issn = {0021-9991},
	doi = {https://doi.org/10.1016/0021-9991(87)90031-3},
	url = {https://www.sciencedirect.com/science/article/pii/0021999187900313},
	author = {Ami Harten and Bjorn Engquist and Stanley Osher and Sukumar R Chakravarthy},
}

@article{onate1993derivation,
	author = {O\~nate, Eugenio and Cervera, Miguel},
	title = {Derivation of thin plate bending elements with one degree of freedom per node: a simple three node triangle},
	journal   = {Eng. Comput.},
	fjournal  = {Engineering Computations},
	volume = {10},
	number = {6},
	pages = {543-561},
	year = {1993},
	month = {06},
	issn = {0264-4401},
	doi = {10.1108/eb023924},
	url = {https://doi.org/10.1108/eb023924},
	eprint = {https://www.emerald.com/ec/article-pdf/10/6/543/702442/eb023924.pdf},
}

@article{hampshire1992three,
	author    = {J. K. Hampshire and B. H. V. Topping and H. C. Chan},
	title     = {Three Node Triangular Bending Elements with One Degree of Freedom per Node},
	journal   = {Eng. Comput.},
	fjournal  = {Engineering Computations},
	year      = {1992},
	volume    = {9},
	number    = {1},
	pages     = {49--62},
	issn      = {0264-4401},
	doi       = {10.1108/eb023848},
	url       = {https://doi.org/10.1108/eb023848},
	publisher = {MCB UP Ltd},
}

@incollection {eymard2000fvm,
	AUTHOR = {Eymard, Robert and Gallou\"et, Thierry and Herbin, Rapha\`ele},
	TITLE = {Finite volume methods},
	BOOKTITLE = {Handbook of numerical analysis, {V}ol. {VII}},
	SERIES = {Handb. Numer. Anal.},
	VOLUME = {VII},
	PAGES = {713--1020},
	PUBLISHER = {North-Holland, Amsterdam},
	YEAR = {2000},
	ISBN = {0-444-50350-1},
	MRCLASS = {65M60 (65N30)},
	MRNUMBER = {1804748},
	MRREVIEWER = {Do\ Young\ Kwak},
	DOI = {10.1016/S1570-8659(00)07005-8},
	URL = {https://doi.org/10.1016/S1570-8659(00)07005-8},
}

@article{godunov1959fvm,
	author = {S. K. Godunov},
	title = {A difference scheme for numerical computation of discontinuous solutions of hydrodynamic equations},
	journal= {Sb. Math.},
	fjournal = {Matematicheskii Sbornik},
	volume = {47},
	year = {1959},
	pages = {271--306}
}

@article{vanLeer1979fvm,
	author = {B. {van Leer}},
	title = {Towards the ultimate conservative difference scheme {V}: A second-order sequel to {G}odunov's method},
	JOURNAL = {J. Comput. Phys.},
	fjournal = {Journal of Computational Physics},
	year = {1979},
	volume = {32},
	number = {1},
	pages = {101--136}
}

@article {hou2024rdg,
	AUTHOR = {Hou, Shijin and Xia, Yinhua},
	TITLE = {Discontinuous {G}alerkin method based on the reduced space for
	the nonlinear convection-diffusion-reaction equation},
	JOURNAL = {J. Sci. Comput.},
	FJOURNAL = {Journal of Scientific Computing},
	VOLUME = {99},
	YEAR = {2024},
	NUMBER = {1},
	PAGES = {Paper No. 19, 25},
	ISSN = {0885-7474,1573-7691},
	MRCLASS = {65M60},
	MRNUMBER = {4715390},
	MRREVIEWER = {Yuan\ Xu},
	DOI = {10.1007/s10915-024-02486-5},
	URL = {https://doi.org/10.1007/s10915-024-02486-5},
}

@book{daniele2012ipdg,
	author    = {Di Pietro, Daniele Antonio and Ern, Alexandre},
	title     = {Mathematical Aspects of Discontinuous {G}alerkin Methods},
	series    = {Math. Appl.},
	volume    = {69},
	publisher = {Springer},
	address   = {Heidelberg},
	year      = {2012},
	isbn      = {978-3-642-22979-4},
	doi       = {10.1007/978-3-642-22980-0}
}

@article {arnold2002ipdg,
	AUTHOR = {Arnold, Douglas N. and Brezzi, Franco and Cockburn, Bernardo
	and Marini, L. Donatella},
	TITLE = {Unified analysis of discontinuous {G}alerkin methods for
	elliptic problems},
	JOURNAL = {SIAM J. Numer. Anal.},
	FJOURNAL = {SIAM Journal on Numerical Analysis},
	VOLUME = {39},
	YEAR = {2001},
	NUMBER = {5},
	PAGES = {1749--1779},
	ISSN = {0036-1429,1095-7170},
	MRCLASS = {65N30},
	MRNUMBER = {1885715},
	DOI = {10.1137/S0036142901384162},
	URL = {https://doi.org/10.1137/S0036142901384162},
}

@incollection{patankar1972fvm,
	title = {{A calculation procedure for heat, mass and momentum transfer in three-dimensional parabolic flows}},
	editor = {Patankar, Suhas V. and Andrew Pollard and Ashok K. Singhal and S. Pratap Vanka},
	booktitle = {Numerical Prediction of Flow, Heat Transfer, Turbulence and Combustion},
	publisher = {Pergamon},
	pages = {54-73},
	year = {1983},
	isbn = {978-0-08-030937-8},
	doi = {https://doi.org/10.1016/B978-0-08-030937-8.50013-1},
	url = {https://www.sciencedirect.com/science/article/pii/B9780080309378500131},
	author = {S. V. Patankar and D. B. Spalding},
}

@article {suli1991poisson,
	AUTHOR = {S\"uli, Endre},
	TITLE = {Convergence of finite volume schemes for {P}oisson's equation
	on nonuniform meshes},
	JOURNAL = {SIAM J. Numer. Anal.},
	FJOURNAL = {SIAM Journal on Numerical Analysis},
	VOLUME = {28},
	YEAR = {1991},
	NUMBER = {5},
	PAGES = {1419--1430},
	ISSN = {0036-1429},
	MRCLASS = {65N06 (65N12 65N15 65N30)},
	MRNUMBER = {1119276},
	MRREVIEWER = {I.\ P.\ Gavrilyuk},
	DOI = {10.1137/0728073},
	URL = {https://doi.org/10.1137/0728073},
}

@book{patankar1980fvm,
	author    = {Patankar, Suhas V.},
	title     = {Numerical Heat Transfer and Fluid Flow},
	publisher = {CRC Press},
	address   = {Boca Raton}, 
	year      = {1980}, 
	note      = {Reprinted edition},
	doi       = {10.1201/9781482234213},
	url       = {https://doi.org/10.1201/9781482234213}
}

@article {mishev1998fvm,
	AUTHOR = {Mishev, Ilya D.},
	TITLE = {Finite volume methods on {V}oronoi meshes},
	JOURNAL = {Numer. Methods Partial Differential Equations},
	FJOURNAL = {Numerical Methods for Partial Differential Equations. An
	International Journal},
	VOLUME = {14},
	YEAR = {1998},
	NUMBER = {2},
	PAGES = {193--212},
	ISSN = {0749-159X,1098-2426},
	MRCLASS = {65N06 (65N12 65N15)},
	MRNUMBER = {1605410},
	MRREVIEWER = {K.\ Najzar},
	DOI = {10.1002/(SICI)1098-2426(199803)14:2<193::AID-NUM4>3.0.CO;2-J},
	URL =
	{https://doi.org/10.1002/(SICI)1098-2426(199803)14:2<193::AID-NUM4>3.0.CO;2-J},
}

@article{chavent1991fvm,
	title = {A unified physical presentation of mixed, mixed-hybrid finite elements and standard finite difference approximations for the determination of velocities in waterflow problems},
	journal ={Adv. Water Resour.},
	fjournal = {Advances in Water Resources},
	volume = {14},
	number = {6},
	pages = {329-348},
	year = {1991},
	issn = {0309-1708},
	doi = {https://doi.org/10.1016/0309-1708(91)90020-O},
	url = {https://www.sciencedirect.com/science/article/pii/030917089190020O},
	author = {G. Chavent and J.E. Roberts}
}

@unpublished{droniou2018gdm,
	TITLE = {{The gradient discretisation method }},
	AUTHOR = {Droniou, J{\'e}r{\^o}me and Eymard, Robert and Gallou{\"e}t, Thierry and Guichard, Cindy and Herbin, Raphaele},
	URL = {https://hal.science/hal-01382358},
	NOTE = {Cf file ''changelog\_hal.pdf''},
	YEAR = {2018},
	MONTH = Mar,
	HAL_ID = {hal-01382358},
	HAL_VERSION = {v6},
}

@article {arbogast1997fvm,
	AUTHOR = {Arbogast, Todd and Wheeler, Mary F. and Yotov, Ivan},
	TITLE = {Mixed finite elements for elliptic problems with tensor
	coefficients as cell-centered finite differences},
	JOURNAL = {SIAM J. Numer. Anal.},
	FJOURNAL = {SIAM Journal on Numerical Analysis},
	VOLUME = {34},
	YEAR = {1997},
	NUMBER = {2},
	PAGES = {828--852},
	ISSN = {0036-1429},
	MRCLASS = {65N30 (65N15 76M10 76S05)},
	MRNUMBER = {1442940},
	MRREVIEWER = {Michael\ J.\ O'Carroll},
	DOI = {10.1137/S0036142994262585},
	URL = {https://doi.org/10.1137/S0036142994262585},
}

@article {ye2004fvm,
	AUTHOR = {Ye, Xiu},
	TITLE = {A new discontinuous finite volume method for elliptic
	problems},
	JOURNAL = {SIAM J. Numer. Anal.},
	FJOURNAL = {SIAM Journal on Numerical Analysis},
	VOLUME = {42},
	YEAR = {2004},
	NUMBER = {3},
	PAGES = {1062--1072},
	ISSN = {0036-1429,1095-7170},
	MRCLASS = {65N30 (35B45 35J55 65N15 76M12)},
	MRNUMBER = {2113675},
	MRREVIEWER = {Qing\ Fang},
	DOI = {10.1137/S0036142902417042},
	URL = {https://doi.org/10.1137/S0036142902417042},
}

@article {ye2007fvm,
	AUTHOR = {Chou, So-Hsiang and Ye, Xiu},
	TITLE = {Unified analysis of finite volume methods for second order
	elliptic problems},
	JOURNAL = {SIAM J. Numer. Anal.},
	FJOURNAL = {SIAM Journal on Numerical Analysis},
	VOLUME = {45},
	YEAR = {2007},
	NUMBER = {4},
	PAGES = {1639--1653},
	ISSN = {0036-1429,1095-7170},
	MRCLASS = {65N06 (35J25 76D07 76M12)},
	MRNUMBER = {2338403},
	MRREVIEWER = {Eugene\ O'Riordan},
	DOI = {10.1137/050643994},
	URL = {https://doi.org/10.1137/050643994},
}

@article {chen2010fvm,
	AUTHOR = {Chen, Long},
	TITLE = {A new class of high order finite volume methods for second
	order elliptic equations},
	JOURNAL = {SIAM J. Numer. Anal.},
	FJOURNAL = {SIAM Journal on Numerical Analysis},
	VOLUME = {47},
	YEAR = {2010},
	NUMBER = {6},
	PAGES = {4021--4043},
	ISSN = {0036-1429,1095-7170},
	MRCLASS = {65N08 (65N12 65N15 65N30)},
	MRNUMBER = {2585177},
	MRREVIEWER = {V\'it\ Dolej\v s\'i},
	DOI = {10.1137/080720164},
	URL = {https://doi.org/10.1137/080720164},
}

@article {bank1987fvm,
	AUTHOR = {Bank, Randolph E. and Rose, Donald J.},
	TITLE = {Some error estimates for the box method},
	JOURNAL = {SIAM J. Numer. Anal.},
	FJOURNAL = {SIAM Journal on Numerical Analysis},
	VOLUME = {24},
	YEAR = {1987},
	NUMBER = {4},
	PAGES = {777--787},
	ISSN = {0036-1429},
	MRCLASS = {65N15 (65N30)},
	MRNUMBER = {899703},
	MRREVIEWER = {Stanis\l aw\ Burys},
	DOI = {10.1137/0724050},
	URL = {https://doi.org/10.1137/0724050},
}

@article {forsyth1988fvm,
	AUTHOR = {Forsyth, Jr., P. A. and Sammon, P. H.},
	TITLE = {Quadratic convergence for cell-centered grids},
	JOURNAL = {Appl. Numer. Math.},
	FJOURNAL = {Applied Numerical Mathematics. An IMACS Journal},
	VOLUME = {4},
	YEAR = {1988},
	NUMBER = {5},
	PAGES = {377--394},
	ISSN = {0168-9274,1873-5460},
	MRCLASS = {65N10},
	MRNUMBER = {948505},
	DOI = {10.1016/0168-9274(88)90016-5},
	URL = {https://doi.org/10.1016/0168-9274(88)90016-5},
}

@article {cao2013fvm,
	AUTHOR = {Cao, Waixiang and Zhang, Zhimin and Zou, Qingsong},
	TITLE = {Superconvergence of any order finite volume schemes for 1{D}
	general elliptic equations},
	JOURNAL = {J. Sci. Comput.},
	FJOURNAL = {Journal of Scientific Computing},
	VOLUME = {56},
	YEAR = {2013},
	NUMBER = {3},
	PAGES = {566--590},
	ISSN = {0885-7474,1573-7691},
	MRCLASS = {65N08 (65L10 65N12)},
	MRNUMBER = {3081662},
	MRREVIEWER = {Erich\ Novak},
	DOI = {10.1007/s10915-013-9691-2},
	URL = {https://doi.org/10.1007/s10915-013-9691-2},
}

@article {cai1991fvem,
	AUTHOR = {Cai, Zhi Qiang},
	TITLE = {On the finite volume element method},
	JOURNAL = {Numer. Math.},
	FJOURNAL = {Numerische Mathematik},
	VOLUME = {58},
	YEAR = {1991},
	NUMBER = {7},
	PAGES = {713--735},
	ISSN = {0029-599X,0945-3245},
	MRCLASS = {65N30 (65N06 65N15)},
	MRNUMBER = {1090257},
	MRREVIEWER = {Lutz\ Angermann},
	DOI = {10.1007/BF01385651},
	URL = {https://doi.org/10.1007/BF01385651},
}

@article {xu2009fvm,
	AUTHOR = {Xu, Jinchao and Zou, Qingsong},
	TITLE = {Analysis of linear and quadratic simplicial finite volume
	methods for elliptic equations},
	JOURNAL = {Numer. Math.},
	FJOURNAL = {Numerische Mathematik},
	VOLUME = {111},
	YEAR = {2009},
	NUMBER = {3},
	PAGES = {469--492},
	ISSN = {0029-599X,0945-3245},
	MRCLASS = {65N06 (65N30)},
	MRNUMBER = {2470148},
	MRREVIEWER = {Srinivasan\ Natesan},
	DOI = {10.1007/s00211-008-0189-z},
	URL = {https://doi.org/10.1007/s00211-008-0189-z},
}

@article {arbogast1998fvm,
	AUTHOR = {Arbogast, Todd and Dawson, Clint N. and Keenan, Philip T. and
	Wheeler, Mary F. and Yotov, Ivan},
	TITLE = {Enhanced cell-centered finite differences for elliptic
	equations on general geometry},
	JOURNAL = {SIAM J. Sci. Comput.},
	FJOURNAL = {SIAM Journal on Scientific Computing},
	VOLUME = {19},
	YEAR = {1998},
	NUMBER = {2},
	PAGES = {404--425},
	ISSN = {1064-8275,1095-7197},
	MRCLASS = {65N30},
	MRNUMBER = {1618879},
	MRREVIEWER = {S.\ Galanis},
	DOI = {10.1137/S1064827594264545},
	URL = {https://doi.org/10.1137/S1064827594264545},
}

@article {cai2003fvm,
	AUTHOR = {Cai, Zhiqiang and Douglas, Jr., Jim and Park, Moongyu},
	TITLE = {Development and analysis of higher order finite volume methods
	over rectangles for elliptic equations},
	JOURNAL = {Adv. Comput. Math.},
	FJOURNAL = {Advances in Computational Mathematics},
	VOLUME = {19},
	YEAR = {2003},
	NUMBER = {1-3},
	PAGES = {3--33},
	ISSN = {1019-7168,1572-9044},
	MRCLASS = {65N06 (65N15)},
	MRNUMBER = {1973457},
	MRREVIEWER = {Miloslav\ Feistauer},
	DOI = {10.1023/A:1022841012296},
	URL = {https://doi.org/10.1023/A:1022841012296},
}
\end{document}